\documentclass[11pt]{amsart}
\usepackage{graphicx,color}

\newtheorem{thm}{Theorem}[section]
\newtheorem{prop}[thm]{Proposition}

\newtheorem{fact}[thm]{Fact}

\newcommand{\pp}[3]{\ensuremath{p_{\overline{#1}#2_{#3}}}}
\newcommand{\qq}[2]{\ensuremath{q_{\overline{#1}#2}}}

\newcommand{\LL}[2]{\ensuremath{L_{\overline{#1}#2}}}

\newcommand{\ZZ}[2]{\ensuremath{Z_{\overline{#1}#2}}}
\newcommand{\bub}[2]{\ensuremath{\mathrm{bub}_{\overline{#1}#2}}}

\newcommand{\CC}[1]{\ensuremath{\mathbf{C}^{#1}}}
\newcommand{\CP}[1]{\ensuremath{\mathbf{CP}^{#1}}}
\newcommand{\PP}[1]{\ensuremath{\mathbf{P}^{#1}}}
\newcommand{\R}[1]{\ensuremath{\mathbf{R}^{#1}}}
\newcommand{\RP}[1]{\ensuremath{\mathbf{RP}^{#1}}}

\newcommand{\GG}[1]{\ensuremath{\Gamma_{#1}}}
\newcommand{\Gb}[1]{\ensuremath{\Gamma_{\overline{#1}}}}
\newcommand{\RR}[2]{\ensuremath{\mathcal{R}_{\overline{#1}#2}}}

\newcommand{\Sym}[1]{\ensuremath{\mathcal{S}_{#1}}}
\newcommand{\Alt}[1]{\ensuremath{\mathcal{A}_{#1}}}
\newcommand{\K}[1]{\ensuremath{\mathcal{K}_{#1}}}
\newcommand{\Kl}[1]{\ensuremath{\mathcal{K}_{#1 \cdot 168}}}
\newcommand{\Z}[1]{\ensuremath{\mathbf{Z}_{#1}}}
\newcommand{\St}[1]{\ensuremath{\mathrm{St}{#1}}}
\newcommand{\Oct}[1]{\ensuremath{\mathcal{O}_{#1}}}
\newcommand{\G}[1]{\ensuremath{\mathcal{G}_{#1}}}
\newcommand{\Ss}[1]{\ensuremath{\Sigma_{#1}}}
\newcommand{\PGL}[2]{\ensuremath{\mathrm{PGL}_{#1}(#2)}}

\newcommand{\pdd}[2]{\frac{\partial{#1}}{\partial{#2}}}

\begin{document}

\title[Solving the heptic in two dimensions]{Solving the heptic in
two dimensions: Geometry and dynamics under Klein's group of order
168}
\author{Scott Crass}
\address{Mathematics Department\\
California State University, Long Beach\\
Long Beach, CA  90840-1001}
 \email{scrass@csulb.edu}
 \keywords{complex dynamics,equivariant map, reflection group,
seventh-degree equation}
\date{\today}

\begin{abstract}

There is a family of seventh-degree polynomials $H$ whose members
possess the symmetries of a simple group of order $168$.  This group
has an elegant action on the complex projective plane.  Developing
some of the action's rich algebraic and geometric properties rewards
us with a special map that also realizes the $168$-fold symmetry.
The map's dynamics provides the main tool in an algorithm that
solves ``heptic" equations in $H$.

\end{abstract}
\maketitle

\section{Introduction and overview}

The alternating group \Alt{7} contains a special subgroup \G{168} of
order 168.  In the 1870s Klein found two linear representations of
\G{168} on \CC{3} that project to an action \K{168} on \CP{2}. As a
subgroup of the symmetric group \Sym{7} whose order is divisible by
seven, \K{168} is a candidate for being the symmetry group of
certain seventh-degree equations.  Following an established course,
we develop an iterative procedure that solves a ``heptic" of this
type. The core of the algorithm is a \K{168}-symmetric map on \CP{2}
whose dynamical output effectively breaks the polynomial's symmetry
and thereby provides for the computation of a root.

The treatment unfolds in stages beginning with a combinatorial
construction of Klein's linear groups \Kl{1} and \Kl{2} that, in
contrast to the discrete group theoretic approach taken by
\cite{klein} and \cite {fricke}, exploits the geometry and
combinatorics of the seven-point projective space and the
three-dimensional octahedral group.  It reproduces many of the
results due to Klein and Fricke while uncovering some new geometric
features. With the action in hand, we develop the groups' respective
systems of invariant differential forms the results of which are
then brought to bear on the rich geometry of \K{168} and the central
question of maps with \K{168} symmetry. Geometrically, a
\K{168}-\emph{equivariant} sends a group orbit to a group orbit. In
algebraic terms, such a map commutes with the action.  The study's
main result is the discovery of a special map $h$ whose strong form
of critical finiteness yields global dynamical properties.  Finally,
$h$ provides for a \emph{reliable} heptic-solving algorithm: for
almost any \K{168}-symmetric heptic and for initial points in an
open dense subset of \CP{2}, the iterative process computes one
(maybe two) of the equation's roots.

\section{Constructing Klein's group}

\subsection{A combinatorial description}

\begin{figure}[h]

\scalebox{1}{\includegraphics{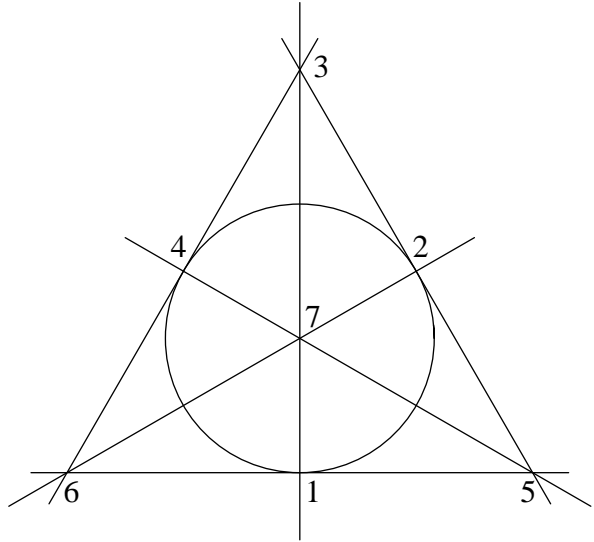}}

\caption{The seven-point projective space $\Z{2}\PP{2}$}

\label{fig:7PtPrjSp}

\end{figure}

The finite projective plane $\Z{2}\PP{2}$ obtained from $\Z{2}^3 =
\Z{2} \times \Z{2} \times \Z{2}$ consists of seven points and seven
lines. Figure~\ref{fig:7PtPrjSp} depicts a well-known representation
of this space where the points appear as $1,\dots,7$ and the lines
as $\overline{124}, \overline{137}, \overline{156}, \overline{235},
\overline{267},\overline{346},\overline{457}$.  (The choice of
labels stems from the arithmetic of \Z{7} in which $1,2,4$ are the
squares and the remaining triples arise by ``translation.") Any pair
of points $a,b$ (lines $L,M$) determines a third point $c$ on
$\overline{abc}$ (line concurrent with $L$ and $M$). Furthermore,
given two points $d,e$ (lines $O,P$), some collineation $\gamma$
that sends $a$ to $d$ and $b$ to $e$ ($L$ to $O$ and $M$ to $P$).
The choice of $a,b$ determines $\gamma(c)$. The image of any one of
the four points (lines) that remain can be any of those four points
(lines).  Choosing the image for one of the four points also
specifies the image of the three remaining points (lines).
Accordingly, the collineation group \PGL{3}{\Z{2}} of $\Z{2}\PP{2}$
is isomorphic to a subgroup of the symmetric group \Sym{7}---indeed,
of the alternating group \Alt{7}---whose order is $168=7\cdot
6\cdot4$.

Selecting a point $a$ (line $\overline{abc}$) divides the lines
(points) into two classes: the three lines through $a$ (points on
$\overline{abc}$) and the four lines that avoid $a$ (points not on
$\overline{abc}$). Under \PGL{3}{\Z{2}}, the combinatorics of these
two sets of lines (points) is isomorphic to that of the octahedron's
three pairs of antipodal vertices and four pairs of antipodal faces
respectively. Applying octahedral motions, the former permute by
\Sym{3} while the latter undergo \Sym{4} behavior. Indeed, \Sym{3}
is the quotient of \Sym{4} by a canonical Klein-$4$ group whose
elements correspond to the rotations by $\pi$ about the three axes
through antipodal vertices. Thus, the stabilizer in \PGL{3}{\Z{2}}
of a point or line is isomorphic to \Sym{4} and permutes the six
complementary points (lines) as the six pairs of antipodal
octahedral edges.

To illustrate, consider the stabilizer \Ss{7} of $7$. The octahedral
configuration---realized by the cube---appears in
Figure~\ref{fig:cubeComb}.  In \Ss{7}, the special Klein-$4$ group
$V_7$ corresponds to the permutations
$$\{(1),(13)(26),(13)(45),(26)(45)\}$$
which pointwise stabilizes the three antipodal face-pairs.  Thus,
$$\Ss{7}/V_7 \simeq \Sym{3}$$
is the group of permutations of the face-pairs. Furthermore, \Ss{7}
contains three conjugate Klein-$4$ subgroups formed, according to
the geometry of $\Z{2}\PP{2}$, by intersections with the stabilizers
of points other than $7$. Furthermore, each of these groups is the
canonical Klein-$4$ subgroup $V_{\overline{ab7}}$ of the stabilizer
of a line $\overline{ab7}$:
\begin{align*}
\{(1),(13)(26),(12)(36),(16)(23)\}&= \Ss{4}\cap \Ss{5}=\Ss{4}\cap
\Ss{7}=\Ss{5}\cap \Ss{7}=
V_{\overline{457}}\\
\{(1),(13)(45),(15)(34),(14)(35)\}&= \Ss{2}\cap
\Ss{6}=\Ss{2}\cap\Ss{7}=\Ss{6}\cap \Ss{7}=
V_{\overline{267}}\\
\{(1),(26)(45),(24)(56),(25)(46)\}&= \Ss{1}\cap
\Ss{3}=\Ss{1}\cap\Ss{7}=\Ss{3}\cap \Ss{7}= V_{\overline{137}}\\.
\end{align*}
(Note that collineations are conflated with their corresponding
permutations.)  For a pair of points, the stabilizer is a dihedral
group $D_4$, but not a subgroup of either point stabilizer. For
instance,
$$D_{17}:=\St{(1,7)}=\{
(1),(26)(45),(24)(56),(25)(46),
(17)(46),(17)(25),(2456)(17),(2654)(17) \}.
$$
The structure among these groups appears in
Figure~\ref{fig:dihedGraph}.
\begin{figure}[h]
\begin{picture}(100,100)(0,0)
\put(50,95){\Ss{\overline{137}}}
\put(10,65){\line(1,1){25}} \put(60,65){\line(0,1){25}}
\put(105,65){\line(-1,1){25}}
\put(0,50){$D_{13}$} \put(50,50){$D_{17}$} \put(100,50){$D_{37}$}
\put(10,40){\line(1,-1){25}} \put(60,40){\line(0,-1){25}}
\put(105,40){\line(-1,-1){25}}
\put(50,0){$V_{\overline{137}}$}
\end{picture}

\caption{Groups associated with the triple of points \{1,3,7\}}

\label{fig:dihedGraph}

\end{figure}
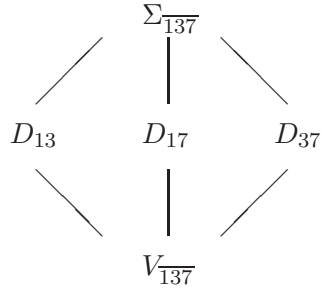
\begin{figure}[ht]

\scalebox{.6}{\includegraphics{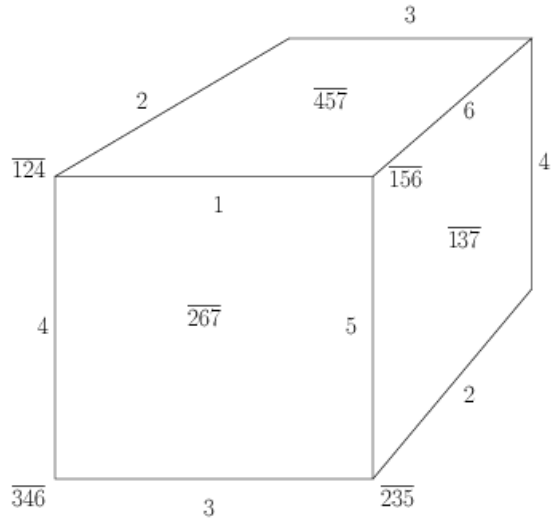}}

\caption{Combinatorics of the cube associated with \Ss{7}}

\label{fig:cubeComb}

\end{figure}

When the stabilizers of a point $a$ and line $\overline{bcd}$ meet,
the result depends on whether or not the point is on the line. If $a
\in \overline{bcd}$, the intersection is a $D_4$. Otherwise, the
result is a $D_3$.  By way of example,
\begin{align*}
\Ss{7} \cap \Ss{\overline{137}}&=D_{13}\\
&=
\{(1),(13)(2465),(26)(45),(13)(2564),(13)(26),(13)(45),(25)(46),(26)(45)\}
\end{align*}
and
$$
\Ss{7} \cap \Ss{\overline{124}}:=D_{\overline{124}7}=
\{(1),(124)(365),(142)(356),(12)(36),(14)(35),(24)(56)\}.
$$

Taking into consideration their intersections, the \Ss{k} account
for $120$ members of \PGL{3}{\Z{2}}. The remaining $48$ elements are
$7$-cycles. To see this, let the $7$-cycle be $(1abcdef)$ and note
that $a$ can be any of the six points other than $1$. Then $b$ is
not on the line containing $1$ and $a$; otherwise, the result would
be the cycle $(1a)$ or $(1ab)$.  So, four choices remain for $b$.
Now, $b$ cannot go to the third point on the line containing $a$ and
$b$---leaving three possibilities for $c$.  However, $1$ and $b$
determine a third point $h$ on $\overline{1bh}$ and $a$ and $h$
determine a third point $g$ on $\overline{agh}$. Accordingly, $c\neq
g$---otherwise, $\overline{1bh}\rightarrow \overline{ach}$ and $h$
would be fixed. Thus, the number of $7$-cycles is $48=6 \cdot 4
\cdot 2$.

Since the \Ss{k} are conjugate, the cycle structure of their
elements determines a conjugacy class for \PGL{3}{\Z{2}}. The
$7$-cycles split into two classes characterized by an element and
its inverse.  Table~\ref{tab:conjCl} summarizes the situation.

\begin{table}[ht]

\begin{tabular}{lccccc}

Cycle structure&
$(ab)(cd)$&$(abcd)(ef)$&$(abc)(def)$&$(abcdefg)$&$(agfedcb)$\\
Number of elements&21&42&56&$24$&$24$

\end{tabular}

\caption{Conjugacy classes for \PGL{3}{\Z{2}}}

\label{tab:conjCl}

\end{table}

\subsection{From combinatorics to geometry}

By situating the cube so that its vertices are $(\pm 1,\pm 1, \pm
1)$, the corresponding octahedral action \Oct{7} on \CP{2} (or
\RP{2}) is expressed by the following elements of
$\mathrm{SO}_3(\CC{})$:
\begin{align*}
&E_1=\begin{pmatrix}
0&0&1\\
0&-1&0\\
1&0&0
\end{pmatrix}
E_2=\begin{pmatrix}
-1&0&0\\
0&0&1\\
0&-1&0
\end{pmatrix}
E_3=\begin{pmatrix}
0&0&-1\\
0&-1&0\\
-1&0&0
\end{pmatrix}
\\
&E_4=\begin{pmatrix}
0&-1&0\\
-1&0&0\\
0&0&-1
\end{pmatrix}
E_5=\begin{pmatrix}
0&1&0\\
1&0&0\\
0&0&-1
\end{pmatrix}
E_6=\begin{pmatrix}
-1&0&0\\
0&0&1\\
0&1&0
\end{pmatrix}\\
&F_{1236}=\begin{pmatrix}
0&-1&0\\
1&0&0\\
0&0&1
\end{pmatrix}
F_{1435}=\begin{pmatrix}
1&0&0\\
0&0&1\\
0&-1&0
\end{pmatrix}
F_{2465}=\begin{pmatrix}
0&0&-1\\
0&1&0\\
1&0&0
\end{pmatrix}\\
&V_{124}=\begin{pmatrix}
0&0&1\\
-1&0&0\\
0&-1&0
\end{pmatrix}
V_{156}=\begin{pmatrix}
0&1&0\\
0&0&1\\
1&0&0
\end{pmatrix}
V_{235}=\begin{pmatrix}
0&1&0\\
0&0&-1\\
-1&0&0
\end{pmatrix}
V_{346}=\begin{pmatrix}
0&-1&0\\
0&0&1\\
-1&0&0
\end{pmatrix}.
\end{align*}
Referring to Figure~\ref{fig:cubeComb}, $E_a$ is an order-$2$
rotation about the $a$ edge-midpoints, $F_{abcd}$ is an order-$4$
rotation about the faces surrounded by $a,b,c,d$, and $V_{abc}$ is
an order-$3$ rotation about the vertex pair attached to $a,b,c$.

By orthogonality, \Oct{7} preserves the quadratic form
$$C_7=x_1^2+x_2^2+x_3^2.$$
Now, the idea is to generate over \Oct{7} an action $\K{168}\simeq
\PGL{3}{\Z{2}}$ on \CP{2}. The procedure determines 1) six quadratic
forms $C_1,\dots, C_6$ other than $C_7$ that \Oct{7} permutes as
octahedral edge-pairs and 2) an order-$4$ transformation $Q \not\in
\Oct{7}$ that permutes the $C_k$---up to a multiplicative
factor---in a way that agrees with the combinatorics of
\PGL{3}{\Z{2}}. First, arbitrarily select $(1326)(57)$ to be the
action of $Q$ on the $C_k$. Then employ undetermined coefficients
for $Q$ and the $C_k$ and enforce the combinatorics of \Oct{7} and
$Q$ on the $C_k$. The result is a system of seven quadratic forms
$C_k$ that realizes the combinatorics of \PGL{3}{\Z{2}} under the
group
$$\K{168}=\bigl< \Oct{7}, Q \bigr>.$$
Specifically, with
$\theta=\tfrac{1+\sqrt{7}\,i}{2}$,
$$
Q=\frac{1}{2}\begin{pmatrix}
1&1&-\theta\\
1&1&\theta\\
1-\theta&\theta-1&0
\end{pmatrix}
$$
and
\begin{align*}
C_1&= \frac{1}{4} \bigl((1-\theta) x_1^2+2 (-2+\theta ) x_2^2+2
(3+\theta ) x_1 x_3+(1-\theta) x_3^2\bigr)\\
C_2&= \frac{1}{4} \bigl(2 (-2+\theta ) x_1^2+(1-\theta) x_2^2-2
(3+\theta ) x_2 x_3+(1-\theta) x_3^2\bigr)\\
C_3&= \frac{1}{4} \bigl((1-\theta) x_1^2+2 (-2+\theta ) x_2^2-2
(3+\theta ) x_1 x_3+(1-\theta) x_3^2\bigr)\\
C_4&= \frac{1}{4} \bigl((1-\theta) x_1^2-2 (3+\theta ) x_1
x_2+(1-\theta) x_2^2+2 (-2+\theta ) x_3^2\bigr)\\
C_5&= \frac{1}{4} \bigl((1-\theta) x_1^2+2 (3+\theta ) x_1
x_2+(1-\theta) x_2^2+2 (-2+\theta ) x_3^2\bigr)\\
C_6&= \frac{1}{4} \bigl(2 (-2+\theta ) x_1^2+(1-\theta) x_2^2+2
(3+\theta ) x_2 x_3+(1-\theta) x_3^2\bigr)\\
C_7&= x_1^2+x_2^2+x_3^2.
\end{align*}
The product
$$
E_1 Q = \frac{1}{2}\begin{pmatrix}
1-\theta&\theta-1&0\\
-1&-1&-\theta\\
1&1&-\theta
\end{pmatrix}
$$
has order $7$ and acts on the $C_k$ by $(1357246)$.  The products
$E_k Q$ and $Q E_k\ (k=1,2,3,6)$ account for the $48$ elements of
order $7$.

It turns out that the matrices given for \Oct{7} and $Q$ permute the
$C_k$ \emph{simply}---no multiplicative character appears.  This
fact implies that
$$\sum_{k=1}^7 C_k \equiv 0.$$
Otherwise, \K{168} would act on the conic $\{\sum_{k=1}^7 C_k =
0\}$---a \CP{1}.  A further consequence concerns linear
representations that project to \K{168}.

\begin{prop}

The projective group \K{168} lifts one-to-one to a group of linear
transformations \Kl{1}.

\end{prop}

\begin{proof}

Lift \Oct{7} and $Q$ to matrices in the literal way and let
$\widehat{\K{}}$ be the group of linear transformations they
generate. Note that $\widehat{\K{}}\subset \mathrm{SL}_3(\CC{})$ so
that the only elements of $\widehat{\K{}}$ that can project to the
identity are $\rho^k I$ where $I$ is the identity matrix, $\rho=e^{2
\pi i/3}$, and $k=0,1,2$.  But, if $k\neq 0$, the action of
$\widehat{\K{}}$ on the $C_k$ gives rise to the character $\rho$ or
$\rho^2$.  Thus, $I$ is the only element of $\widehat{\K{}}$ that
projects to the identity and the projection is one-to-one.

\end{proof}

We can extend \Kl{1} to a double cover of \K{168}:
$$\Kl{2}=\Kl{1} \cup -\Kl{1}.$$
For each of the $21$ involutions $Z$ in \Kl{1}, $-Z$ is conjugate to
$$\begin{pmatrix}
1&0&0\\
0&1&0\\
0&0&-1
\end{pmatrix}$$
and so pointwise fixes a plane---as well as its defining linear
form---and setwise fixes a line in \CC{3}. In \CP{2}, the projective
version of $Z$ fixes a line and point. Since the $21$
transformations $-Z$ generate \Kl{2}, this group possesses the rich
structure of a \emph{complex reflection group}. (See \cite{st}.)

The combinatorial agreement between \PGL{3}{\Z{2}} and \K{168} gives
rise to a second system of seven quadratic forms and conics: the
sums of the triples $C_a+C_b+C_c$ where $\overline{abc}$ is a line
in $\Z{2}\PP{2}$.  The resulting forms are:
\begin{align*}
C_{\overline{124}}&= \frac{1}{2} \bigl(-x_1^2-x_2^2-(3+\theta ) x_1
\bigl(x_2-x_3\bigr)-(3+\theta ) x_2 x_3-x_3^2\bigr)\\
C_{\overline{137}}&= \frac{1}{2} \bigl((3-\theta) x_1^2+2
(-1+\theta ) x_2^2+(3-\theta) x_3^2\bigr)\\
C_{\overline{156}}&= \frac{1}{2} \bigl(-x_1^2-x_2^2+(3+\theta ) x_2
x_3-x_3^2+(3+\theta ) x_1 \bigl(x_2+x_3\bigr)\bigr)\\
C_{\overline{235}}&= \frac{1}{2} \bigl(-x_1^2-x_2^2+(3+\theta ) x_1
\bigl(x_2-x_3\bigr)-(3+\theta ) x_2 x_3-x_3^2\bigr)\\
C_{\overline{267}}&= \frac{1}{2} \bigl(2 (-1+\theta )
x_1^2+(3-\theta) \bigl(x_2^2+x_3^2\bigr)\bigr)\\
C_{\overline{346}}&= \frac{1}{2} \bigl(-x_1^2-x_2^2+(3+\theta ) x_2
x_3-x_3^2-(3+\theta ) x_1 \bigl(x_2+x_3\bigr)\bigr)\\
C_{\overline{457}}&= \frac{1}{2} \bigl((3-\theta) x_1^2+(3-\theta)
x_2^2+2 (-1+\theta ) x_3^2\bigr).
\end{align*}
Attention now turns to the elegant algebraic and geometric
structures connected to these forms and their associated conics.

\section{Algebra of invariants}

In the search for special maps with \K{168} symmetry, the algebra
and geometry of \K{168} invariant forms are illuminating. A basic
tool of classical invariant theory is the \emph{Molien series}: a
formal power series associated with a linear representation of a
finite group whose coefficients express the dimension of the space
of polynomial invariants in the associated degree.  By extension,
there is a two-variable Molien series for the invariant differential
forms under the group's action on the exterior polynomial algebra
$\Lambda(\CC{n})$. Of course, the space of polynomials is identified
with the $0$-forms.  Given an invariant $p$-form $\alpha$ of degree
$\ell$, the exterior derivative $d \alpha$ is an invariant
$(p+1)$-form of degree $\ell-1$. Also, for an invariant $q$-form
$\beta$ of degree $m$, $\alpha \wedge \beta$ is an invariant
$(p+q)$-form of degree $\ell+m$.

Molien's Theorem provides a generating function for the the series.

\begin{thm}

Let $G$ be a representation on \CC{n} of a finite group.  Then the
extended Molien series is
$$
M_G(\Lambda(\CC{n}))=\sum_{k=0}^n \biggl( \sum_{\ell=1}^\infty
\dim{(\Lambda_{k,\ell}^G)}\, t^\ell \biggr) s^k
 = \frac{1}{|G|}\sum_{A\in G} \frac{\det{(I+s\,A)}}{\det{(I-t\,A)}}
 = \frac{1}{|G|}\sum_{A\in \mathcal{T}} |\mathcal{C}_A|
 \frac{\det{(I+s\,A)}}{\det{(I-t\,A)}}
$$
where $\Lambda_{k,\ell}^G$ is the space of $G$-invariant $k$-forms
with polynomial degree $\ell$, $\mathcal{T}$ is a transversal for
the conjugacy classes of $G$, and $\mathcal{C}_A$ is the conjugacy
class of $A$.

\end{thm}
Application of the theorem to the lifts of \K{168} to \CC{3} yields
interesting and, eventually, useful results.

\begin{fact}
\begin{align*}
M_{\Kl{2}}(\Lambda(\CC{3}))=&\
\frac{(1+st^3)(1+ st^5)(1+st^{13})}{(1-t^4)(1-t^6)(1-t^{14})} \\
=&\ (1 + t^4 + t^6 + t^8 + t^{10} + 2\,t^{12} + 2\,t^{14} +
2\,t^{16} + 3\,t^{18} + 3\,t^{20}+\dots ) s^0
\\
&+  (t^3 + t^5 + t^7 + 2\,t^9 + 2\,t^{11} + 3\,t^{13} + 3\,t^{15} +
5\,t^{17} + 5\,t^{19} + 6\,t^{21}+\dots) s^1   \\
&+ ( t^8 + t^{12} + t^{14} + 2\,t^{16} + 2\,t^{18} + 3\,t^{20} +
4\,t^{22} + 4\,t^{24} + 5\,t^{26}+\dots) s^2
\nonumber\\
&+ t^{21} (1 + t^4 + t^6 + t^8 + t^{10} + 2\,t^{12} + 2\,t^{14} +
2\,t^{16} + 3\,t^{18} + 3\,t^{20} +\dots) s^3.
\end{align*}
\end{fact}

The infinite nature of these series is due to the process of
promoting a $k$-form $\alpha$ of degree-$\ell$ to a $k$-form of
degree $\ell+m$ under multiplication by a $0$-form of degree $m$. In
particular, the effect of promoting \Kl{1}-invariant $k$-forms by
\K{2\cdot 168}-invariant $0$-forms disappears upon dividing the
series for the former by that for the latter.  The result is a
finite sum.

\begin{fact}
\begin{align*}
\frac{M_{\Kl{1}}(\Lambda(\CC{3}))}
{M_{\Kl{2}}(\Lambda_0(\CC{3}))}=&\ (1 + t^{21}) s^0 \\
&+ (t^3 + t^5 + t^{10} + t^{12} + t^{13} + t^{20}) s^1 \\
&+( t + t^8 + t^9 + t^{11} + t^{16} + t^{18}) s^2  \\
&+(1 + t^{21}) s^3.
\end{align*}
\end{fact}
Rearranging the factors yields the generating function for the
simple lift of \K{168}:
\begin{align*}
M_{\Kl{1}}(\Lambda(\CC{3})) =&\
\frac{1}{(1-t^4)(1-t^6)(1-t^{14})} \bigl((1 + t^{21}) (s^0+s^3) \\
&+ (t^3 + t^5 + t^{10} + t^{12} + t^{13} + t^{20}) s^1 +(t + t^8 +
t^9
+ t^{11} + t^{16} + t^{18}) s^2\bigr)\\
=&\ (1 + t^4 + t^6 + t^8 + t^{10} + 2\,t^{12} + 2\,t^{14} +
2\,t^{16} + 3\,t^{18} + 3\,t^{20} + t^{21} +3\,t^{22} +\dots)
s^0 \\
&+ (t^3 + t^5 + t^7 + 2\,t^9 + t^{10} + 2\,t^{11} + t^{12} +
3\,t^{13} + t^{14} + 3\,t^{15} + 2\,t^{16}+\dots ) s^1  \\
&+(t + t^5 + t^7 + t^8 + 2\,t^9 + 2\,t^{11} + t^{12} + 3\,t^{13} +
t^{14} + 4\,t^{15} + 2\,t^{16}+\dots ) s^2  \\
&+ (1 + t^4 + t^6 + t^8 + t^{10} + 2\,t^{12} + 2\,t^{14} + 2\,t^{16}
+ 3\,t^{18} + 3\,t^{20} + t^{21} +\dots) s^3.
\end{align*}

As for interpreting these results, consider first the series for the
reflection group \Kl{2}.  Manifesting a result of reflection group
theory, its generating function indicates that the invariant
$0$-forms are, under the three-dimensional action \Kl{2}, polynomial
combinations of three basic forms---call them $F$, $\Phi$, and
$\Psi$---whose respective degrees are $4$, $6$, and $14$.  The
invariant $1$-forms are generated as a module over the $0$-forms by
the three exterior derivatives $dF$, $d\Phi$, and $d\Psi$ whose
degrees are $3$, $5$, and $13$. In a similar way, the three products
$dF \wedge d\Phi$, $dF \wedge d\Phi$, and $d\Phi \wedge d\Psi$ of
degrees $8$, $16$, and $18$ generate the module of $2$-form
invariants.  Finally, the invariant $3$-forms are the multiples of
$dF \wedge d\Phi \wedge d\Psi$. Here another general reflection
group property is realized: when expressed in coordinates as
$$dF \wedge d\Phi \wedge d\Psi=X(x) dx_1 \wedge dx_2 \wedge dx_3,$$
the degree-$21$ expression $X$ is, up to a  constant multiple, the
product of the $21$ linear forms fixed by the involutions that
generate \Kl{2}.  This fact accounts for the factor $t^{21}$ in the
Molien series. Of course, $X=\det{J(F,\Phi,\Psi)}$ where $J(f,g,h)$
is the Jacobian matrix of $f$, $g$, and $h$.  Note that $X$ is not a
\K{2\cdot 168} invariant, but rather a \emph{relative} invariant
where
$$
X\circ T =
\begin{cases}
X& T\in \Kl{1}\\
-X& T\not\in \Kl{1}
\end{cases}.
$$

From the generating function for \Kl{1} we glean that $F$, $\Phi$,
$\Psi$, and $X$ generate the ring of invariants $\CC{}[x]^{\Kl{1}}$.
By the relative invariance of $X$, $\CC{}[x]^{\K{1\cdot 168}}$
admits a relation
$$X^2=P(F,\Phi,\Psi)$$
for some polynomial $P$ that will be specified below.  In the case
of $1$-forms, three new entries $\zeta$, $\eta$, and $\theta$ appear
in degrees $10$, $12$, and $20$.  Evidently, $\theta=dX$. However,
$\zeta$ and $\eta$ are not exterior derivatives so that $d\zeta$ and
$d\eta$ do not vanish identically and account for the $2$-forms of
degrees $9$ and $11$.  The remaining item is the degree-$1$
``identity $2$-form" expressed by
$$
x_1 dx_2 \wedge dx_3 + x_2 dx_3 \wedge dx_1 + x_3 dx_1 \wedge dx_2.
$$
Notice the numerological duality between $1$-forms and $2$-forms:
each $1$-form has a companion $2$-form such that the degrees sum to
$21$; to wit, $3+18$, $5+16$, etc.  For an invariant theoretic
account of this phenomenon see \cite{shepler}.

\clearpage
\section{Geometry}

\subsection{Configurations}

\subsubsection{Octahedral conics}

Each quadratic form $C_k$ and $C_{\overline{abc}}$ defines a conic
$$
\GG{k}=\{C_k=0\} \quad \text{and} \quad
\Gb{abc}=\{C_{\overline{abc}}=0\}
$$
as a \CP{1} endowed with the structure of an octahedron.  Due to
their connection to \PGL{3}{\Z{2}}, call the \GG{k} ``point-conics"
and the \GG{\overline{abc}} ``line-conics."  In parallel to the
treatment of the combinatorics under $\mathrm{PGL}_3(\Z{2})$, two
types of intersection can occur \emph{between} these systems of
conics. When $k\in\{a,b,c\}$ (for definiteness, let $k=b$), $\GG{b}
\cap \Gb{abc}$ is a two-point set $\{\pp{ac}{b}{1,2}\}$ that
determines a $D_4$-stable line \LL{ac}{b}.  In fact, \LL{ac}{b} is
pointwise fixed by an involution \ZZ{ac}{b} that lifts to one of the
$21$ generating reflections in \K{2\cdot 168}. The six points on
\GG{b} given by intersections of this type are octahedral vertices;
overall, they belong to an orbit of $42$ points and determine the
orbit of $21$ mirrors for the reflection group \K{2\cdot 168}. (For
ease of reference, call such a point or line a ``$42$-point" or
``$21$-line" and employ this style of terminology to orbits of any
size.)  The dihedral structure here appears when two octahedra are
``joined" along an axis through a pair of antipodal vertices.

In case $k\not\in \{a,b,c\}$, the intersection consists of two
points \pp{abc}{k}{1,2} on a line \LL{abc}{k} with a $D_3$
stabilizer. Intersections of this kind give, on a single conic,  the
eight face-centers of the octahedron and, over all conics, orbits of
$56$-points and $28$-lines.  Here, two ``intersecting" octahedra
meet along an axis through an antipodal pair of face-centers at an
angular separation of $180^\circ$.  This $D_3$
structure---illustrated in Figure~\ref{fig:bubConfig}---will play a
heuristic role in determining a symmetry that exchanges the two
system of conics.

\begin{figure}[ht]

\scalebox{.85}{\includegraphics{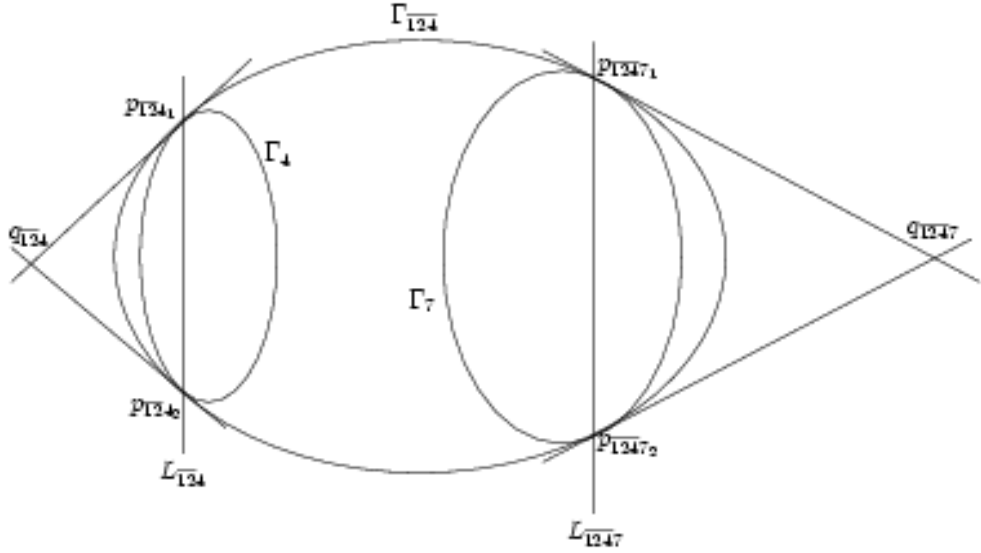}}

\caption{Intersecting conics}

\label{fig:conicInt}

\end{figure}

Furthermore, the lines tangent to \GG{b} at a pair of $42$-points
$\pp{ac}{b}{1,2}$ intersect in a $21$-point \qq{ac}{b}---the pole of
\LL{ac}{b} with respect to \GG{b}---fixed by \ZZ{ac}{b}.  A similar
situation produces a triangle of two $56$-points and a $28$-point.
Figure~\ref{fig:conicInt} illustrates the configuration.

As an example, Table~\ref{tab:conicLines} indicates how special
lines fall on \GG{7}.  The $12$ edge-midpoints on \GG{7} belong to a
\K{168} orbit of $84$-points. Where do the lines tangent to \GG{7}
at a pair of antipodal edge-midpoints meet?  Facts bear out the
notation's suggestion: the $21$-point \qq{37}{1} lies on each of the
lines tangent to \GG{3} and \GG{7} at their respective intersections
with \LL{37}{1}.  The twelve remaining $21$-lines hit \GG{7} in an
octahedral orbit of 24 points.  Such an orbit is typically generic,
but this one has special significance.  Are they the vertices of a
special truncation of the octahedron?

\begin{table}[ht]

\begin{tabular}{cl}

\hline

Line&Octahedral axis\\

\hline

\LL{13}{7}\ \LL{26}{7}\ \LL{45}{7} &Antipodal vertices\\

\LL{124}{7}\ \LL{156}{7}\ \LL{235}{7}\ \LL{346}{7}&
Antipodal faces\\

$\underbrace{
 \underbrace{\LL{37}{1}\ \LL{17}{3}}\
 \underbrace{\LL{27}{6}\ \LL{67}{2}}\
 \underbrace{\LL{47}{5}\ \LL{57}{4}}}_{
 \text{``great circles" of edges}}$
&Antipodal edges\\

\hline\\

\end{tabular}

\caption{Special lines through \GG{7}}

\label{tab:conicLines}

\end{table}

\subsubsection{Real projective octahedra} An octahedron defines a
natural antipodal structure on a sphere. Relative to an octahedral
conic, a point $X$ and its antipode $\widehat{X}$ determine a polar
point $P(X)$ in \CP{2}. The collection
$$\RR{}{k}=\{P(X)\ |\ X\in \GG{k}\}$$
is an octahedral \RP{2} with $3$ vertices, $4$ faces, and $6$ edges.
For \RR{}{7}, the preceding discussion shows that

\begin{itemize}

\item  the vertices are $21$-points: $\qq{13}{7}, \qq{26}{7},
\qq{45}{7}$

\item the face-centers are $28$-points: $\qq{124}{7},
\qq{156}{7}, \qq{235}{7}, \qq{346}{7}$

\item the edge-midpoints are $21$-points:
$\qq{37}{1}, \qq{17}{3}, \qq{67}{2}, \qq{27}{6}, \qq{57}{4},
\qq{47}{5}$.

\end{itemize}

\begin{figure}[ht]

\scalebox{1.1}{\includegraphics{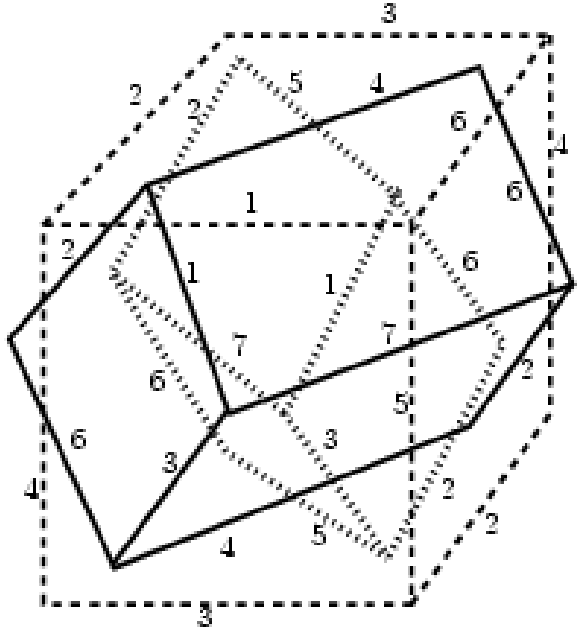}}

\caption{Combinatorial arrangement for three octahedral conics}

\label{fig:3Cubes}

\end{figure}

Combinatorial considerations showed that two octahedral groups in
the \emph{same} system \Oct{a} and \Oct{b} intersect in one of the
three conjugate Klein-$4$ subgroups of each octahedral group. We can
realize this structure with polyhedra by superimposing two cubes
each one representing an octahedral group by means of the
combinatorics of its edge-pairs.  That is, the edge-pairs of the
``$a$-cube" and ``$b$-cube" are labeled by $\{1,\dots,7\}-\{a\}$ and
$\{1,\dots,7\}-\{b\}$ respectively.  The triple $\overline{abc}$
determines a third $c$-cube. Arrange the three cubes in the
following way.  (See Figure~\ref{fig:3Cubes}.)

\begin{itemize}

\item  The $a$-edge on the $b$-cube ($c$-cube) cuts diagonally across
the face of the $b$-cube ($c$-cube) whose edge-labels are
$\{1,\dots,7\}-\{a,b,c\}$.

\item The $b$-edge on the $a$-cube ($c$-cube) cuts diagonally across
the face of the $a$-cube ($c$-cube) whose edge-labels are
$\{1,\dots,7\}-\{a,b,c\}$.

\item The $c$-edge on the $a$-cube ($b$-cube) cuts diagonally across
the face of the $a$-cube ($b$-cube) whose edge-labels are
$\{1,\dots,7\}-\{a,b,c\}$.

\end{itemize}

\begin{figure}[ht]

\scalebox{1}{\includegraphics{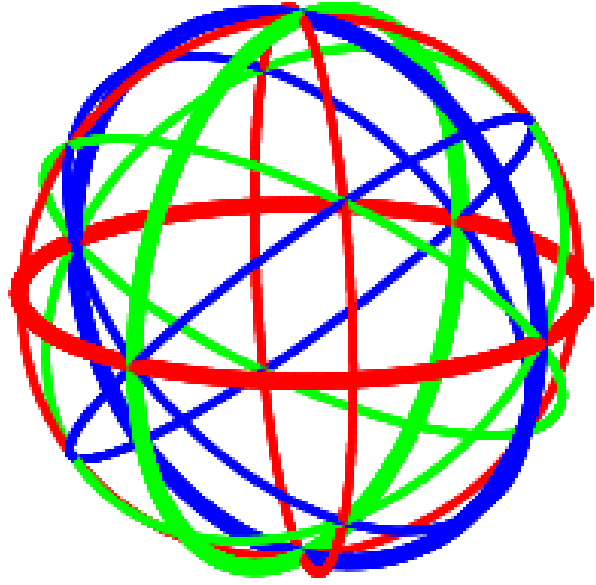}}

\caption{Three octahedra forming a fourth}

\label{fig:3Oct}

\end{figure}

To illustrate, let $\{a,b,c\}=\{4,5,7\}$.  The three groups \Oct{4},
\Oct{5}, \Oct{7} meet in the Klein-$4$ group whose action on edges
of each cube is
$$V_{\overline{457}}=\{(1),(12)(36),(16)(23),(13)(26)\}.$$
Applied to the configuration of cubes on a sphere, the edges labeled
$4$, $5$, and $7$---two pairs from each cube---yield three mutually
orthogonal great circles that form an octahedron. By dualizing the
three cubes to octahedra, we find the geometric manifestation of the
combinatorial phenomenon depicted in Figure~\ref{fig:dihedGraph}. A
model of this structure appears in Figure~\ref{fig:3Oct}. Each
octahedron---red, blue, green---contributes one (thick) great circle
to produce a fourth octahedron whose symmetry group is
$\Oct{\overline{457}}=\St{(\GG{4} \cup \GG{5} \cup \GG{7})}$.

\subsubsection{Reflection lines}

The notation for a $21$-line employs a pairing of a point and
line-conic. This mode of reference supplies answers to several
questions. How does the collection of $21$-lines intersect a given
$21$-line \LL{ac}{b}?  Which special points belong to a $21$-line?
How does the reflection \ZZ{ac}{b} act on the octahedral conics?

Being $D_4$-stable, a $21$-point \qq{ac}{b} is fixed by four
reflections other than \ZZ{ac}{b}. Hence, four $21$-lines contain
\qq{ac}{b} with an intersection multiplicity of $6=\binom{4}{2}$.
Similarly, a $D_3$-stable $28$-point \qq{abc}{d} lies on three such
lines.  These concurrences account for all intersections of
$21$-lines.  In numerical terms, the total number of intersections
with multiplicity is
$$\binom{21}{2}=21 \binom{4}{2} + 28
\binom{3}{1}.$$ Accordingly, \LL{ac}{b} contains four $21$-points
and four $28$-points.  The notation readily tells us which ones.
First, consider the $7\times 7$ array of Table~\ref{tab:21Array}
that distinguishes $21$ entries in an obvious way.


\begin{table}[ht]

$\begin{array}{c|c|c|c|c|c|c|c|}
&1&2&3&4&5&6&7\\
\hline
&&&&&&&\\
 \overline{124}&
  \overline{24}{1}&\overline{14}{2}&&\overline{12}{4}&&&\\
\hline
&&&&&&&\\
 \overline{137}&
  \overline{37}{1}&&\overline{17}{3}&&&&\overline{13}{7}\\
\hline
&&&&&&&\\
 \overline{156}&
  \overline{56}{1}&&&&\overline{16}{5}&\overline{15}{6}&\\
\hline
&&&&&&&\\
 \overline{235}&
  &\overline{35}{2}&\overline{25}{3}&&\overline{23}{5}&&\\
\hline
&&&&&&&\\
 \overline{267}&
  &\overline{67}{2}&&&&\overline{27}{6}&\overline{26}{7}\\
\hline
&&&&&&&\\
 \overline{346}&
  &&\overline{46}{3}&\overline{36}{4}&&\overline{34}{6}&\\
\hline
&&&&&&&\\
 \overline{457}&
  &&&\overline{57}{4}&\overline{47}{5}&&\overline{45}{7}\\
\hline
\end{array}$

\caption{A combinatorial array}

\label{tab:21Array}

\end{table}

For the $\overline{25}{3}$ entry, its row and column each contain
two other entries indicating which $21$-points belong to \LL{25}{3}
(or, by a \K{168} duality, which $21$-lines contain \qq{25}{3}):
$$
\qq{17}{3},\qq{46}{3},\qq{35}{2},\qq{23}{5} \in \LL{25}{3} \qquad
 \LL{17}{3} \cap \LL{46}{3} \cap \LL{35}{2} \cap \LL{25}{3} =
 \{\qq{25}{3}\}.
$$
The rows and columns containing these four entries intersect in four
empty cells that correspond to the $28$-points on \LL{25}{3},
namely, \qq{137}{2}, \qq{137}{5}, \qq{346}{2}, and \qq{346}{5}.
Finally, the row and column of the $\overline{124}{7}$ cell each
have three entries.  For each such row, there's exactly one column
with a shared entry. These three entries give the $21$-lines on a
$28$-point:
$$\LL{37}{1} \cap \LL{67}{2} \cap \LL{57}{4} = \{\qq{124}{7}\}.$$
In the dual situation, the three respective poles \qq{37}{1},
\qq{67}{2}, and \qq{57}{4} of these lines relative to the associated
pair of conics lie on \LL{124}{7}.

As for the action on conics by \ZZ{25}{3}, the column entries
$\overline{17}{3}$ and $\overline{46}{3}$ indicate that the
permutation of point-conics and line-conics undergo $(17)(46)$ and
$(\overline{124}\ \overline{267})(\overline{156}\ \overline{457}).$

\subsection{Antiholomorphic symmetry and special coordinates}

\begin{figure}[ht]

\scalebox{.75}{\includegraphics{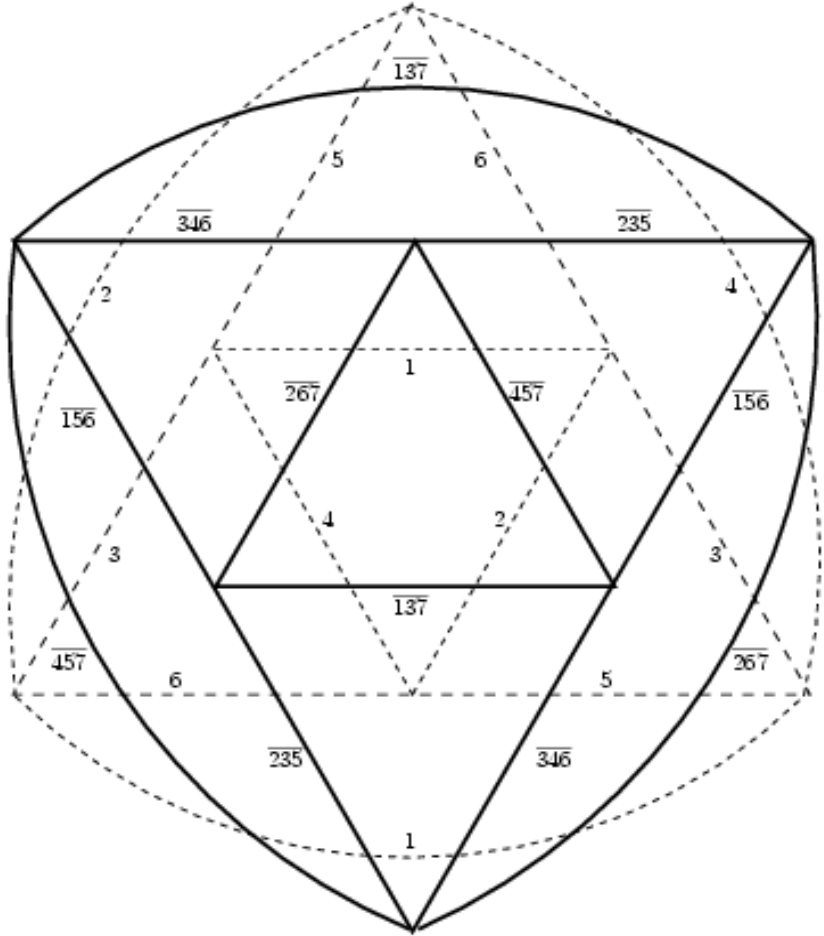}}

\caption{Heuristic for antiholomorphic exchange of conic systems}

\label{fig:bubConfig}

\end{figure}

As an analogue to the icosahedral and Valentiner groups, \K{168}
extends to a group of order $336$ by an antiholomorphic involution
that exchanges the two systems of octahedral conics. (See
\cite{sextic}, \S~2C.)  Adopting the approach taken in the
Valentiner case, we look for a heuristic that suggests a permutation
between point and line-conics.  (\cite{sextic}, \S~2D)  Recall the
device of regarding the intersection of a point and complementary
line-conic as a superposition of octahedra that ``meet" at an
antipodal pair of face-centers.  Figure~\ref{fig:bubConfig} shows
such an intersection for \GG{7} (solid line segments) and \Gb{124}
(dashed segments) as planar nets whose edge-labels corresponding to
the respective conics.  Treating this image as two octahedra viewed
from above the ``intersecting" faces---one pair at the north pole
and the other at the south, a reflection in the equatorial plane
exchanges the octahedra and acts on the edges as follows:
$$
1\leftrightarrow \overline{137} \quad
 2\leftrightarrow \overline{267} \quad
  4\leftrightarrow \overline{457} \quad
3\leftrightarrow \overline{156} \quad
 6\leftrightarrow \overline{235} \quad
  5\leftrightarrow \overline{346}.
$$
This pairing of point and line-conics also emerges through
examination of the $7\times 7$ array.

An antiholomorphic involution that exchanges conics in the suggested
way---of course, $7\leftrightarrow \overline{124}$---would fix the
$28$-points \qq{124}{7}, \qq{156}{3}, \qq{235}{6}, and \qq{346}{5}.
For historical reasons, such a map is called \bub{124}{7}. A
candidate for \bub{124}{7} is the conjugation map
$$
[y_1,y_2,y_3] \longrightarrow
[\overline{y_1},\overline{y_2},\overline{y_3}],
$$
a condition that requires the specified points to be expressible by
real coordinates.  Accordingly, making a change to $y$ coordinates
in which these points are symmetrically arranged to be
$$
[0,0,1],[1,0,1],
 \biggl[-\frac{1}{2},\frac{\sqrt{3}}{2},1\biggr],
 \biggl[-\frac{1}{2},-\frac{\sqrt{3}}{2},1\biggr]
$$
respectively results in the prescribed exchange of point and
line-conics under conjugation of coordinates.  To give the simplest
expression:
$$
C_7(y)=y_1^2+y_2^2+\frac{7}{8}(1-3\,\theta)y_3^2 \qquad \text{and}
\qquad C_{\overline{124}}(y)=\overline{C_7(\overline{y})}.
$$

\begin{figure}[ht]

\scalebox{.85}{\includegraphics{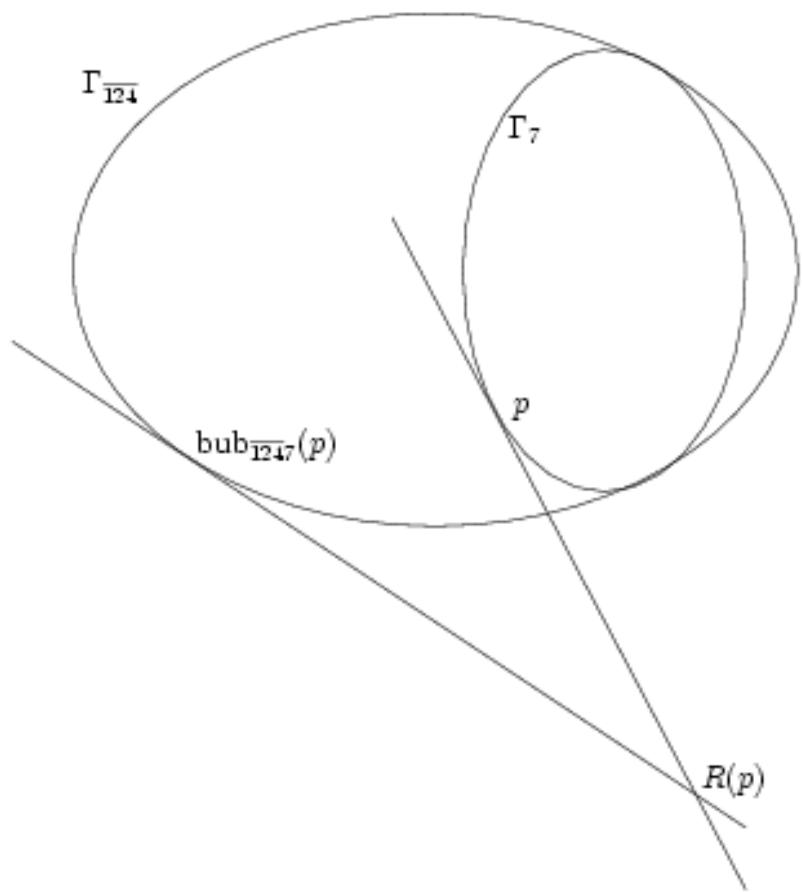}}

\caption{Construction for a reflection \RP{2}}

\label{fig:bubGeom}

\end{figure}

The above-mentioned degree-two extension of \K{168} appears as
$$\overline{\K{}}_{336}=\bigl<\K{168},\bub{124}{7}\bigr>.$$
Each of the $28$ point and line-conic pairs $(\GG{k},\Gb{abc})$,
$k\not\in \{a,b,c\}$, is canonically associated with an
antiholomorphic involution \bub{abc}{k} the collection of which
generates $\overline{\K{}}_{336}$. Moreover, each \bub{abc}{k}
pointwise fixes an \RP{2}.  For \bub{124}{7}, the ``mirror" is
$$
 \RR{124}{7}=\{[y_1,y_2,y_3]\ |\ y_1,y_2,y_3\in \R{}\}
$$
which, as Figure~\ref{fig:bubGeom} shows, is geometrically
determined by the pair of conics $(\Gb{124},\GG{7})$.  Let $p\in
\GG{7}$ and $T_p\GG{7}$ be the line tangent to \GG{7} at $p$. Then
the point
$$R(p)=T_p\GG{7} \cap T_{\bub{124}{7}(p)}\Gb{124}$$
is fixed by \bub{124}{7}.  As $p$ ranges over \GG{7}, $R(p)$ covers
the $D_3$-symmetric \RR{124}{7}.  Since there is a second point on
\GG{7} whose tangent line passes through $R(p)$, the map $R:\GG{7}
\rightarrow \RR{124}{7}$ is, as it should be, two-to-one.

Upcoming discussions will further explore the geometry of the \RP{2}
mirrors and will benefit from use of $y$ coordinates.

\newpage

\subsection{Klein's quartic and other basic invariants}

The Molien series reveals the presence of a unique lowest-degree
\Kl{2} invariant in degree four.  Since \Kl{2} permutes the
octahedral quadratic forms $C_k$, symmetric functions in the $C_k$
are \Kl{2} invariant.  In particular, $\sum C_k^2$ gives a non-zero
polynomial; when expressed in $y$ and a common numerical factor is
dropped, the result is
$$
 F=
  3\,(y_1^4 + y_2^4) - 49\,y_3^4
  -14\,y_1^3 y_3 +6\,y_1^2 y_2^2 - 21\,(y_1^2 +y_2^2)y_3^2
  +42\,y_1 y_2^2 y_3.
$$
The unique sixth-degree invariant is given by $\sum C_k^3$.
Alternatively, the determinant of the Hessian matrix $H(F)$ of $F$
produces this form:
\begin{align*}
\Phi=&\ \alpha_{\Phi}\,\det{H(F)}\\
=&\ 26 y_1^6+210 y_3 y_1^5+330 y_2^2 y_1^4+315 y_3^2 y_1^4+980 y_3^3
y_1^3-420 y_2^2 y_3 y_1^3-90 y_2^4 y_1^2-3675 y_3^4 y_1^2 \\
&+630 y_2^2 y_3^2 y_1^2-2940 y_2^2 y_3^3 y_1-630 y_2^4 y_3 y_1+54
y_2^6+686 y_3^6-3675 y_2^2 y_3^4+315 y_2^4 y_3^2.
\end{align*}
Here, $\alpha_{\Phi}$ removes the highest common factor of the
coefficients.

A third and final generator $\Psi$ for the \Kl{2} invariants results
from applying the classical tool of the ``bordered Hessian."
Removing a common integer factor $\tfrac{1}{\alpha_{\Psi}}$,
\begin{align*} \Psi=&\ \alpha_{\Psi}\,B(F,\Phi)=
\alpha_{\Psi}\,\det{\begin{pmatrix}
&&&\pdd{\Phi}{y_1}\\[5pt]
&H(F)&&\pdd{\Phi}{y_2}\\[5pt]
&&&\pdd{\Phi}{y_3}\\[5pt]
\pdd{\Phi}{y_1}&\pdd{\Phi}{y_2}&\pdd{\Phi}{y_3}&1
\end{pmatrix}}\\
=&\ 4388 y_1^{14}+62916 y_3 y_1^{13}-343756 y_2^2 y_1^{12}+378574
y_3^2 y_1^{12}+853384 y_3^3 y_1^{11}-297528 y_2^2 y_3 y_1^{11}
\\
&-267036 y_2^4 y_1^{10}-4618152 y_3^4 y_1^{10}-3304560 y_2^2 y_3^2
y_1^{10}-26185306 y_3^5 y_1^9+15510460 y_2^2 y_3^3 y_1^9 \\
&+391020 y_2^4 y_3 y_1^9-1432620 y_2^6 y_1^8-48569829 y_3^6
y_1^8-78795675 y_2^2 y_3^4 y_1^8+14116410 y_2^4 y_3^2 y_1^8
\\
&+70186032 y_3^7 y_1^7+43131564 y_2^2 y_3^5 y_1^7-44205840 y_2^4
y_3^3 y_1^7+529200 y_2^6 y_3 y_1^7-605556 y_2^8 y_1^6 \\
&-90656958 y_3^8 y_1^6-56457114 y_2^2 y_3^6 y_1^6-120454740 y_2^4
y_3^4 y_1^6-27222048 y_2^6 y_3^2 y_1^6+156473170 y_3^9 y_1^5
\\
&-70186032 y_2^2 y_3^7 y_1^5-1037232 y_2^4 y_3^5 y_1^5-34747272
y_2^6 y_3^3 y_1^5-1666980 y_2^8 y_3 y_1^5+927612 y_2^{10} y_1^4
\\
&+237533331 y_3^{10} y_1^4-897745905 y_2^2 y_3^8 y_1^4-245478240
y_2^4 y_3^6 y_1^4-33802650 y_2^6 y_3^4 y_1^4 \\
&+23694930 y_2^8 y_3^2 y_1^4-448007392 y_3^{11} y_1^3-312946340
y_2^2 y_3^9 y_1^3-350930160 y_2^4 y_3^7 y_1^3 \\
&+300537972 y_2^6 y_3^5 y_1^3+16669800 y_2^8 y_3^3 \ y_1^3-1492344
y_2^{10} y_3 y_1^3-6804 y_2^{12} y_1^2+130119794 y_3^{12} y_1^2
\\
&+475066662 y_2^2 y_3^{10} y_1^2+145212480 y_2^4 y_3^8
y_1^2-270847206 y_2^6 y_3^6 y_1^2+1666980 y_2^8 y_3^4 y_1^2
\\
&+952560 y_2^{10} y_3^2 y_1^2+1344022176 y_2^2 y_3^{11}
y_1-469419510 y_2^4 y_3^9 y_1-210558096 y_2^6 y_3^7 y_1 \\
&+64178730 y_2^8 y_3^5 y_1-7445844 y_2^{10} y_3^3 y_1-47628 y_2^{12}
y_3 y_1-26244 y_2^{14}+3294172 y_3^{14} \\
&+130119794 y_2^2 y_3^{12}+237533331 y_2^4 y_3^{10}-160187517 y_2^6
y_3^8-33256251 y_2^8 y_3^6 \\
&-10807587 y_2^{10} y_3^4+23814 y_2^{12} y_3^2.
\end{align*}

According to an earlier discussion, the extra invariant $X$ for
\Kl{1} is given by the Jacobian of $F$, $\Phi$, and $\Psi$ as well
as by the product of the $21$ linear forms that define the
reflection lines.  With a normalizing constant $\alpha_X$,
\begin{align*}
X=&\ \alpha_X\,\det{J(F,\Phi,\Psi)}\\
=&\ y_2 (3 y_1^2-y_2^2) (7 y_1^2-14 y_3 y_1+3 y_2^2+7 \ y_3^2) (4
y_1^2+7 y_3 y_1+7 y_3^2) (2 y_1^2-7 y_3 y_1+14 y_3^2) \\
&(16 y_1^4+56 y_3 y_1^3+36 y_2^2 y_1^2+105 y_3^2 y_1^2+98 y_3^3
y_1+36 y_2^4+49 y_3^4-63 y_2^2 y_3^2) \\
&(4 y_1^4-28 y_3 y_1^3-24 y_2^2 y_1^2+105 y_3^2 y_1^2-196 y_3^3
y_1+84 y_2^2 y_3 y_1+36 y_2^4+196 y_3^4+21 y_2^2 y_3^2) \\
& (y_1^4+14 y_3 y_1^3-6 y_2^2 y_1^2+105 y_3^2 y_1^2+392 y_3^3 y_1-42
y_2^2 y_3 y_1+9 y_2^4+784 y_3^4+21 y_2^2 y_3^2).
\end{align*}
Computing the relation among the basic \Kl{1} invariants,
\small
\begin{align*} X^2=&\ \frac{1}{823543}(23328 F^9 \Phi + 9288
F^6 \Phi^3 + 938 F^3 \Phi^5 + \Phi^7 - 3888 F^7 \Psi - 612 F^4
\Phi^2 \Psi \\
&+ 21 F \Phi^4 \Psi - 66 F^2 \Phi \Psi - \Psi^3).
\end{align*}
\normalsize

\begin{figure}[ht]

\scalebox{1}{\includegraphics{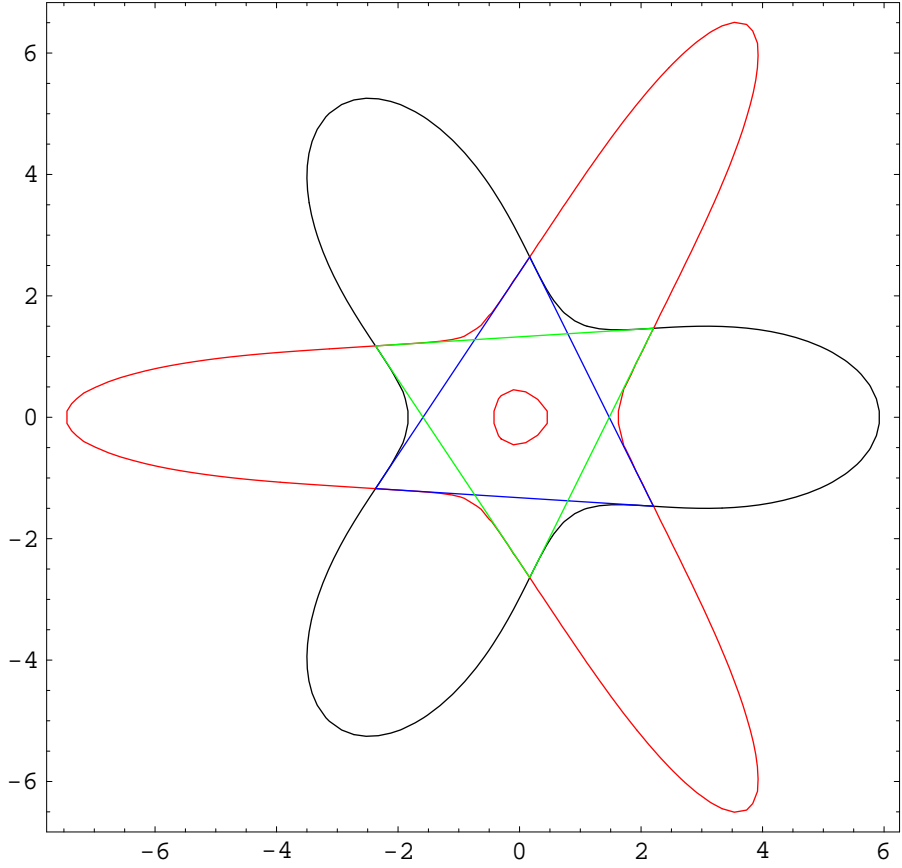}}

\caption{Klein's quartic in a \bub{}{}-fixed \RP{2}}

\label{fig:quartic}

\end{figure}

The curve $\mathcal{F}=\{F=0\}$ is the celebrated Klein quartic. To
touch on the rich geometry of this surface, its inflection points
lie on $\mathcal{G}=\{\Phi=0\}$ and form a $24$-point \K{168} orbit.
A $24$-point is fixed by an order-seven transformation so that the
orbit decomposes into eight triples---the three points fixed by each
\Z{7} subgroup. Reproducing one of Klein's pictures,
Figure~\ref{fig:quartic} displays $\mathcal{F}$ (black) and
$\mathcal{G}$ (red) on the affine plane $\{y_3\neq 0\}$ in
\RR{124}{7} in the symmetrical $y$ coordinates. In this space,
$\mathcal{F}\cap\mathcal{G}$ consists of two triples of $24$-points
each of which forms a triangle of lines that are tangent to one of
the curves at a $24$-point.  The curves' three axes of symmetry are
\RP{1} intersections of the three $21$-lines associated with
\RR{124}{7}. (See \cite{8fold} for a broad treatment of Klein's
quartic.)

\subsection{Catalog of special orbits}

The following table summarizes some of the \K{168} geometry.

\begin{table}[ht]

\begin{tabular}{cccl}

\hline

Descriptor&Size&Stabilizer&Geometric data\\
\hline
 \qq{ac}{b}&$21$&$D_4$&
  On four $21$-lines\\
  &&&Vertices of an octahedral
  \RP{2}\\
  &&&Edge midpoints of two octahedral \RP{2}s\\
 &$24$&\Z{7}&$\mathcal{F}\cap\mathcal{G}$\\
 \qq{abc}{k}&$28$&$D_3$&
 On three $21$-lines\\
 &&&Face-centers of an octahedral \RP{2}\\
\pp{ac}{b}{1,2}&$42$&\Z{4}&
 Vertices of an octahedral conic (\CP{1})\\
\pp{abc}{k}{1,2}&$56$&\Z{3}&Face-centers of an octahedral conic\\
&&&$\mathcal{F}\cap\{\Psi=0\}$\\
&$84$&\Z{2}&Generic points on $21$-lines\\

\hline

\end{tabular}

\end{table}

\section{Maps with \K{168} symmetry}

\subsection{Equivariant basics}

Here's a useful fact: under an action \G{} on \CC{3},
\G{}-equivariant maps are in one-to-one correspondence with
\G{}-invariant $2$-forms. (\cite{sextic}, \S\,3A)  Thus, the Molien
series can help uncover maps that generate the module of
\Kl{1}-equivariant maps over \Kl{1}-invariant polynomials.  There
are six basic equivariants in degrees $1$, $8$, $9$, $11$, $16$, and
$18$.  By a previous discussion on invariant $2$-forms, the even
degree maps are given by
$$
\psi=\nabla F \times \nabla \Phi \qquad
 \phi=\nabla F \times \nabla \Psi \qquad
  f=\nabla \Phi \times \nabla \Psi
$$
where $\nabla$ indicates a formal gradient operator and $\times$ is
the cross product.  For example, expressed in coordinates
$$
 dF \wedge d\Phi = (\nabla F \times \nabla \Phi)^T \begin{pmatrix}
 dy_2\wedge dy_3 \\
 dy_3\wedge dy_1 \\
 dy_1\wedge dy_2
 \end{pmatrix}.
$$
Moreover, these three maps generate the \Kl{2} equivariants. As for
the odd-degree maps, the only obvious case is the identity. Whether
there is an exterior algebraic procedure that produces the maps
$g_9$ and $g_{11}$ in degrees $9$ and $11$ is unknown. Here, our
concern is computing these maps.

Since $X\cdot g_9$ is a degree-$30$ \Kl{2} equivariant, it can be
expressed by a combination of suitably promoted versions of $\psi$,
$\phi$, and $f$.  That is, for some $a_k$,
$$X\cdot g_9=
 (a_1\,F^4 \Phi + a_2\,F \Phi^3 + a_3\,F^2 \Psi) \psi +
  (a_4\,F^2 \Phi + a_5\,\Psi) \phi +
   (a_6\,F^3 + a_7\,\Phi^2) f.
$$
By solving a system of linear equations, we can choose the $a_k$ so
that the ``$30$-map" vanishes at $\{X=0\}$ and so, is divisible by
$X$. Dividing by $X$ and a common integer factor produces 
\tiny
\begin{align*} g_9=& [27 y_1^9+336 y_3 y_1^8-4428
y_2^2 y_1^7+308 y_3^2 y_1^7-14700 y_3^3 y_1^6-10920 y_2^2 y_3 
y_1^6-13446 y_2^4 y_1^5-151851 y_3^4 y_1^5+18564 y_2^2 y_3^2 
y_1^5-149548 y_3^5 y_1^4\\
&-192864 y_2^2 y_3^3 y_1^4+33600 y_2^4 y_3 y_1^4-5436 y_2^6 y_1^3 
-8232 y_3^6 y_1^3-44394 y_2^2 y_3^4 y_1^3-116676 y_2^4 y_3^2 
y_1^3+605052 y_3^7 y_1^2-329280 y_2^2 y_3^5 y_1^2\\
&-104076 y_2^4 y_3^3 y_1^2-11592 y_2^6 y_3 y_1^2+3555 y_2^8 
y_1-261709 y_3^8 y_1-8232 y_2^2
y_3^6 y_1-238287 y_2^4 y_3^4 y_1+37548 y_2^6 y_3^2 y_1-605052 y_2^2 
y_3^7\\
&+259308 y_2^4 y_3^5+74088 y_2^6 y_3^3, \\
&y_2 (4509 y_1^8+4956 y_3 y_1^7-1620 y_2^2 y_1^6+39382 y_3^2
y_1^6-207564 y_3^3 y_1^5-22596 y_2^2 y_3 y_1^5-7722 y_2^4 y_1^4 \\
&+22197 y_3^4 y_1^4-111174 y_2^2 y_3^2 y_1^4-628376 y_3^5
y_1^3-177576 y_2^2 y_3^3 y_1^3+24612 y_2^4 y_3 y_1^3-612 y_2^6 y_1^2 
+8232 y_3^6 y_1^2+476574 y_2^2 y_3^4 y_1^2\\
&+24066 y_2^4 y_3^2 y_1^2+1210104 y_3^7 y_1-189336 y_2^2 y_3^5 
y_1+29988 y_2^4 y_3^3 y_1-4284 y_2^6 y_3 y_1+981 y_2^8+261709 
y_3^8+8232 y_2^2 y_3^6\\
&+108633 y_2^4 y_3^4+2142 y_2^6 y_3^2),\\
&224 y_1^9+2547 y_3 y_1^8-2016 y_2^2 y_1^7+9156 y_3^2 y_1^7+7952
y_3^3 y_1^6+19260 y_2^2 y_3 y_1^6+6048 y_2^4 y_1^5-57918 y_3^4 y_1^5
-9156 y_2^2 y_3^2 y_1^5+150675 y_3^5 y_1^4\\
&-43176 y_2^2 y_3^3 y_1^4+18306 y_2^4 y_3 y_1^4-6048 y_2^6 
y_1^3-48706 y_3^6 y_1^3 +115836 y_2^2 y_3^4 y_1^3-45780 y_2^4 y_3^2 
y_1^3-51450 y_3^7 y_1^2+301350 y_2^2 y_3^5 y_1^2\\
&+68544 y_2^4 y_3^3 y_1^2+5148 y_2^6 y_3 y_1^2 +146118 y_2^2 y_3^6 
y_1+173754 y_2^4 y_3^4 y_1-27468 y_2^6 y_3^2 y_1+7203 y_3^9-51450 
y_2^2 y_3^7+150675 y_2^4 y_3^5\\
&+504 y_2^6 y_3^3+3555 y_2^8 y_3].
\end{align*}
\normalsize
By a similar process in degree $32=21+11$, 
\tiny
\begin{align*}
g_{11}=& [-1334 y_1^{11}-2534 y_3 y_1^{10}-63034 y_2^2 y_1^9+47355
y_3^2 y_1^9+542136 y_3^3 y_1^8+45864 y_2^2 y_3 y_1^8-88716 y_2^4
y_1^7+1184232 y_3^4 y_1^7-551460 y_2^2 y_3^2 y_1^7\\
&-312816 y_3^5 y_1^6+2446080 y_2^2 y_3^3 y_1^6-308196 y_2^4 y_3 
y_1^6-114084 y_2^6 y_1^5-4519368 y_3^6 y_1^5-6856668 y_2^2 y_3^4 
y_1^5+735714 y_2^4 y_3^2 y_1^5\\
&+432180 y_3^7 y_1^4+625632 y_2^2 y_3^5 y_1^4-446880 y_2^4 y_3^3 
y_1^4+634536 y_2^6 y_3 y_1^4+18306 y_2^8 y_1^3-11503191 y_3^8 
y_1^3+4272408 y_2^2 y_3^6 y_1^3\\
&+4515840 y_2^4 y_3^4 y_1^3+1066716 y_2^6 y_3^2 y_1^3+5176556 y_3^9 
y_1^2+18670176 y_2^2 y_3^7 y_1^2+938448 y_2^4 y_3^5 y_1^2-268128 
y_2^6 y_3^3 y_1^2-162918 y_2^8 y_3 y_1^2\\
&+5022 y_2^{10} y_1+9311078 y_3^{10} y_1-11503191 y_2^2 y_3^8 
y_1-8956416 y_2^4 y_3^6 y_1+2019780 y_2^6 y_3^4 y_1-267813 y_2^8 
y_3^2 y_1-5176556 y_2^2 y_3^9\\
&-6655572 y_2^4 y_3^7-814968 y_2^8 
y_3^3 ,\\
&y_2 (-17294 y_1^{10}-60746 y_3 y_1^9-36546 y_2^2 y_1^8+309309 y_3^2
y_1^8-1631112 y_3^3 y_1^7+356328 y_2^2 y_3 y_1^7-77628 y_2^4 
y_1^6-1056930 y_3^4 y_1^6\\
&-900732 y_2^2 y_3^2 y_1^6-312816 y_3^5 y_1^5-5794152 y_2^2 y_3^3 
y_1^5-579852 y_2^4 y_3 y_1^5+52236 y_2^6 y_1^4+2136204 y_3^6 
y_1^4-1627290 y_2^2 y_3^4 y_1^4\\
&-486738 y_2^4 y_3^2 y_1^4-19534536 y_3^7 y_1^3+625632 y_2^2 y_3^5 
y_1^3-2348472 y_2^4 y_3^30 y_1^3+168840 y_2^6 y_3 y_1^3+1674 y_2^8 
y_1^2-11503191 y_3^8 y_1^2\\
&-17912832 y_2^2 y_3^6 y_1^2+9745218 y_2^4 y_3^4 y_1^2+717444 y_2^6 
y_3^2 y_1^2-10353112 y_3^9 y_1+5359032 y_2^2 y_3^7 y_1+938448 y_2^4 
y_3^5 y_1-1083096 y_2^6 y_3^3 y_1\\
&+11718 y_2^8 y_3 y_1-8586 y_2^{10}+9311078 y_3^{10}-11503191 y_2^2 
y_3^8-2300844 
y_2^4 y_3^6-221382 y_2^6 y_3^4-5859 y_2^8 y_3^2),\\
&2464 y_1^{11}+8670 y_3 y_1^{10}-19712 y_2^2 y_1^9+26250 y_3^2 
y_1^9-196119 y_3^3 y_1^8-79878 y_2^2 y_3 y_1^8+44352 y_2^4
y_1^7-89376 y_3^4 y_1^7-201096 y_2^2 y_3^2 y_1^7\\
&-1953924 y_3^5 y_1^6-276444 y_2^2 y_3^3 y_1^6-77604 y_2^4 y_3 
y_1^6+1422078 y_3^6 y_1^5+89376 y_2^2 y_3^4 y_1^5+579852 y_2^4 y_3^2 
y_1^5+1825446 y_3^7 y_1^4\\
&-6454476 y_2^2 y_3^5 y_1^4-1007370 y_2^4 y_3^3 y_1^4+114084 y_2^6 
y_3 y_1^4-66528 y_2^8 y_1^3-5978490 y_3^8 y_1^3-2844156 y_2^2 y_3^6 
y_1^3+446880 y_2^4 y_3^4 y_1^3\\
&-634536 y_2^6 y_3^2 y_1^3+4477865 y_3^9 y_1^2+3650892 y_2^2 y_3^7 
y_1^2-5466636 y_2^4 y_3^5 y_1^2-1066716 y_2^6 y_3^3 y_1^2+98118 
y_2^8 y_3 y_1^2+17935470 y_2^2 y_3^8 y_1\\
&-4266234 y_2^4 y_3^6 y_1+268128 y_2^6 y_3^4 y_1-11718 y_2^8 y_3^2 
y_1+100842 y_3^{11}+4477865 y_2^2 y_3^9+1825446 y_2^4 y_3^7-2019780 
y_2^6 y_3^5\\
&-139671 y_2^8 y_3^3-5022 y_2^{10} y_3].
\end{align*}
\normalsize

\subsection{A special map}

For a group action \G, an equivariant $g$, and a point $a$ fixed by
$T\in\G{}$,
$$Tg(a)=g(Ta)=g(a).$$
Accordingly, a point fixed by \ZZ{ac}{b} has an equivariant image
that's also fixed by \ZZ{ac}{b}. Thus, a $21$-line \LL{ac}{b} either
maps to itself or collapses to its companion point \qq{ac}{b}.
Under a holomorphic map, the only possible outcome is the line's
preservation.

On numerological grounds, the family of $15$-maps offers a special
opportunity.  Since the critical set of a $15$-map has degree
$42=3(15-1)$, there could be a map whose critical set is the
$21$-lines with multiplicity two.  Such an outcome would follow a
pattern found in the \Sym{n} actions on \CP{n-2}. (See
\cite{critFin}.)

The \Kl{1} equivariant $15$-maps are given by
$$
g_{15}=(a_1\,F^2 \Phi + a_2\,\Psi)\cdot \epsilon + a_3\,\Phi\cdot
g_9 + a_4\,F\cdot g_{11}
$$
where $\epsilon$ denotes the identity map. Taking the $a_k$ so that
$\det{J(g_{15})}|_{\{X=0\}}=0$ yields 
\tiny
\begin{align*}
h=&\ \biggl(\frac{827}{14}\,F^2\Phi +
\frac{99}{49}\,\Psi\biggr)\cdot \epsilon -
 \frac{11}{9882516}\,\Phi\cdot g_9 +
 \frac{3}{1647086}\,F\cdot g_{11} \\
=&\ [27 (816 y_1^{15}+4200 y_3 y_1^{14}-18960 y_2^2 y_1^{13}+1680
y_3^2 y_1^{13}+97020 y_3^3 y_1^{12}-68040 y_2^2 y_3 y_1^{12}-4640 
y_2^4 y_1^{11}-332955 y_3^4 y_1^{11}\\
&-150780 y_2^2 y_3^2 y_1^{11}-950796 y_3^5 y_1^{10}+776160 y_2^2 
y_3^3 y_1^{10}+80080 y_2^4 y_3 y_1^{10}-96800 y_2^6 y_1^9-2716560 
y_3^6 y_1^9-7858620 y_2^2 y_3^4 y_1^9\\
&+412720 y_2^4 y_3^2 y_1^9+2376990 y_3^7 y_1^8+679140 y_2^2 y_3^5 
y_1^8-2651880 y_2^4 y_3^3 y_1^8+314160 y_2^6 y_3 y_1^8-34320 y_2^8 
y_1^7-1188495 y_3^8 y_1^7\\
&-339570 y_2^2 y_3^6 y_1^7-9588810 y_2^4 y_3^4 y_1^7-1321320 y_2^6 
y_3^2 y_1^7+11092620 y_3^9 y_1^6-1584660 y_2^4 y_3^5 y_1^6-1811040 
y_2^6 y_3^3 y_1^6-120120 y_2^8 y_3 y_1^6\\
&+77616 y_2^{10} y_1^5+7764834 y_3^{10} y_1^5-83194650 y_2^2 y_3^8 
y_1^5-5546310 y_2^4 y_3^6 y_1^5-1584660 y_2^6 y_3^4 y_1^5+776160 
y_2^8 y_3^2 y_1^5-38824170 y_3^{11} y_1^4\\
&-11092620 y_2^2 y_3^9 y_1^4-7923300 y_2^4 y_3^7 y_1^4+7923300 y_2^6 
y_3^5 y_1^4+743820 y_2^8 y_3^3 y_1^4-138600 y_2^{10} y_3 
y_1^4+52689945 y_2^4 y_3^8 y_1^3\\
&-11884950 y_2^6 y_3^6 y_1^3+1754445 y_2^8 y_3^4 y_1^3-97020 
y_2^{10} y_3^2 y_1^3+77648340 y_2^2 y_3^{11} y_1^2-25882780 y_2^4 
y_3^9 y_1^2-6338640 y_2^6 y_3^7 y_1^2\\
&+2716560 y_2^8 y_3^5 y_1^2-776160 y_2^{10} y_3^3 y_1^2+12941390 
y_2^4 y_3^{10} y_1-11092620 y_2^6 y_3^8 y_1-3961650 y_2^8 y_3^6 
y_1-1358280 y_2^{10} y_3^4 y_1\\
&+12941390 y_2^4 y_3^{11}-3697540 y_2^6 y_3^9-792330 y_2^8 
y_3^7+271656 y_2^{10}
y_3^5),\\
&-27 y_2^3 (24640 y_1^{12}+61600 y_3 y_1^{11}-8624 y_2^2
y_1^{10}+204820 y_3^2 y_1^{10}-1196580 y_3^3 y_1^9-154000 y_2^2 y_3
y_1^9+66000 y_2^4 y_1^8+7413945 y_3^4 y_1^8\\
&-582120 y_2^2 y_3^2 y_1^8-2942940 y_3^5 y_1^7+2587200 y_2^2 y_3^3 
y_1^7-147840 y_2^4 y_3 y_1^7+26400 y_2^6 y_1^6+3961650 y_3^6 
y_1^6+13809180 y_2^2 y_3^4 y_1^6\\
&+1265880 y_2^4 y_3^2 y_1^6+6338640 y_3^7 y_1^5-2852388 y_2^2 y_3^5 
y_1^5+3039960 y_2^4 y_3^3 y_1^5+110880 y_2^6 y_3 y_1^5-76320 y_2^8 
y_1^4+85967805 y_3^8 y_1^4\\
&+18223590 y_2^2 y_3^6 y_1^4+1730190 y_2^4 y_3^4 y_1^4-665280 y_2^6 
y_3^2 y_1^4+29580320 y_3^9 y_1^3+12677280 y_2^2 y_3^7 y_1^3-10639860 
y_2^4 y_3^5 y_1^3-258720 y_2^6 y_3^3 y_1^3\\
&+171360 y_2^8 y_3 y_1^3-2160 y_2^{10} y_1^2-25882780 y_3^{10} 
y_1^2-49916790 y_2^2 y_3^8 y_1^2+15507030 y_2^4 y_3^6 y_1^2-1067220 
y_2^6 y_3^4 y_1^2+49140 y_2^8 y_3^2 y_1^2\\
&-103531120 y_3^{11} y_1+29580320 y_2^2 y_3^9 y_1+6338640 y_2^4 
y_3^7 y_1-2037420 y_2^6 y_3^5 y_1+485100 y_2^8 y_3^3 y_1-15120 
y_2^{10} y_3 y_1+1296 y_2^{12}\\
&-5176556 y_2^2 y_3^{10}+9904125 y_2^4 y_3^8+1245090 y_2^6 
y_3^6+963585 y_2^8 y_3^4+7560
y_2^{10} y_3^2),\\
&-704 y_1^{15}+36960 y_3 y_1^{14}+10560 y_2^2 y_1^{13}+221760 y_3^2
y_1^{13}+560560 y_3^3 y_1^{12}-406560 y_2^2 y_3 y_1^{12}-63360 y_2^4
y_1^{11}+2651880 y_3^4 y_1^{11}\\
&-1552320 y_2^2 y_3^2 y_1^{11} -7844067 y_3^5 y_1^{10}-4107180 y_2^2 
y_3^3 y_1^{10}+1552320 y_2^4 y_3 y_1^{10}+190080 y_2^6 
y_1^9-23185085 y_3^6 y_1^9+18078060 y_2^2 y_3^4 y_1^9\\
&+2217600 y_2^4 y_3^2 y_1^9-68083785 y_3^7 y_1^8-199667160 y_2^2 
y_3^5 y_1^8+12806640 y_2^4 y_3^3 y_1^8-1995840 y_2^6 y_3 
y_1^8-285120 y_2^8 y_1^7+121226490 y_3^8 y_1^7\\
&-6791400 y_2^2 y_3^6 y_1^7-57047760 y_2^4 y_3^4 y_1^7+3991680 y_2^6 
y_3^2 y_1^7-187914265 y_3^9 y_1^6+24109470 y_2^2 y_3^7 
y_1^6-292369770 y_2^4 y_3^5 y_1^6\\
&-26195400 y_2^6 y_3^3 y_1^6-997920 y_2^8 y_3 y_1^6+171072 y_2^{10} 
y_1^5+202440315 y_3^{10} y_1^5-121226490 y_2^2 y_3^8 y_1^5+164012310 
y_2^4 y_3^6 y_1^5\\
&-61122600 y_2^6 y_3^4 y_1^5-5987520 y_2^8 y_3^2 y_1^5+150506685 
y_3^{11} y_1^4-2402344560 y_2^2 y_3^9 y_1^4-309687840 y_2^4 y_3^7 
y_1^4-42785820 y_2^6 y_3^5 y_1^4\\
&+27941760 y_2^8 y_3^3 y_1^4+2993760 y_2^{10} y_3 y_1^4-767071480 
y_3^{12} y_1^3-404880630 y_2^2 y_3^{10} y_1^3-606132450 y_2^4 y_3^8 
y_1^3+171143280 y_2^6 y_3^6 y_1^3\\
&-1746360 y_2^8 y_3^4 y_1^3-5987520 y_2^{10} y_3^2 y_1^3-141178800 
y_3^{13} y_1^2+301013370 y_2^2 y_3^{11} y_1^2+661991715 y_2^4 y_3^9 
y_1^2-437026590 y_2^6 y_3^7 y_1^2\\
&+32089365 y_2^8 y_3^5 y_1^2+2619540 y_2^{10} y_3^3 y_1^2+2301214440 
y_2^2 y_3^{12} y_1-607320945 y_2^4 y_3^{10} y_1-363679470 y_2^6 
y_3^8 y_1+71819055 y_2^8 y_3^6 y_1\\
&-13097700 y_2^{10} y_3^4 y_1+105413504 y_3^{15}-141178800 y_2^2 
y_3^{13} +150506685 y_2^4 y_3^{11}-392203350 y_2^6 y_3^9-35145495 
y_2^8 y_3^7-25671492 y_2^{10} y_3^5] .
\end{align*} \normalsize
The condition for vanishing on $\{X=0\}$ implies that $\det{J(h)}= G
X$ where $G$ is a degree-$21$ \Kl{1} invariant.  Since $X$ is the
only such invariant, $\det{J(h)}$ is divisible by $X^2$ as well.

By preserving $21$-lines, $h$ acquires a strong form of critical
finiteness; to wit, the restriction $h|_{\LL{ac}{b}}$ is also
critically finite with critical points where the other $21$-lines
meet \LL{ac}{b}, namely, the $21$ and $28$-points. Arguments adapted
from \cite{critFin} for maps of this type establish two significant
results.

\begin{thm}

Under \K{168}, $h$ is the unique holomorphic equivariant that's
doubly critical exactly at the $21$-lines.

\end{thm}

Furthermore, $h$ possesses global dynamical properties that qualify
it for the leading role in a \emph{reliable} heptic-solving
algorithm

\begin{thm}

The $21$ and $28$-points are superattracing fixed points whose
basins exhaust the Fatou set $F_h$ and are dense in \CP{2}.

\end{thm}

\begin{figure}[ht]

\resizebox{\textwidth}{!}{\includegraphics{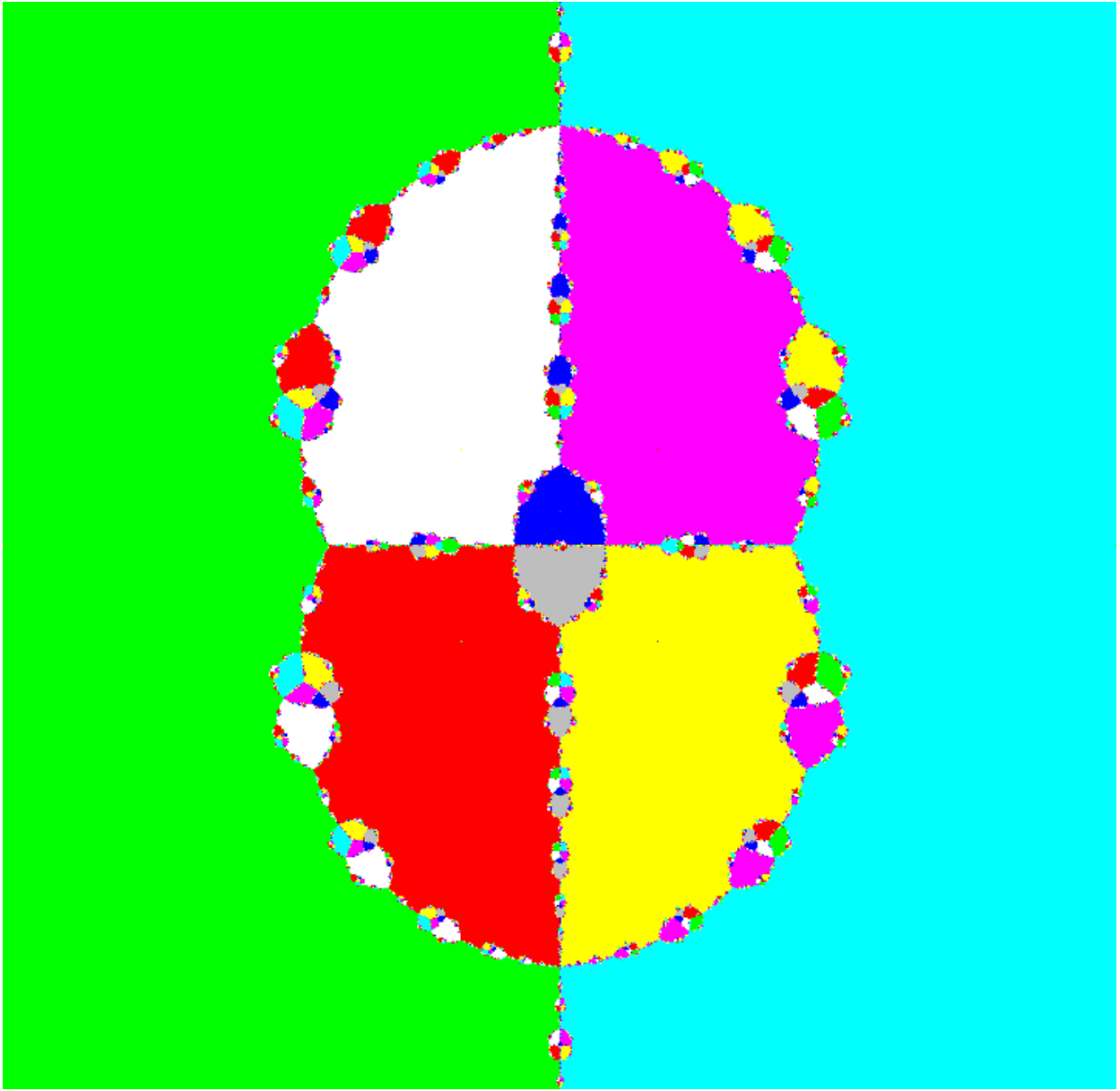}}

\caption{Basins of attraction on a $21$-line}

\label{fig:hL21}

\end{figure}

The geometric character of $h$ appears in its basins of attraction.
Figure~\ref{fig:hL21} shows the basin plot for
$\widetilde{h}=h|_{\LL{37}{1}}$ where we see the line as \CC{}.
Symmetrical coordinates are chosen so that the $21$ and $28$-points
are
\begin{gather*}
 \qq{17}{3}=e^{\pi i/4} \quad
 \qq{24}{1}=e^{3\,\pi i/4} \quad
 \qq{13}{7}=e^{5\,\pi i/4} \quad
 \qq{56}{1}=e^{7\,\pi i/4}\\
 \qq{124}{7}=\frac{\sqrt{2}}{3-\sqrt{7}} \quad
 \qq{156}{7}=\frac{\sqrt{2}\,i}{3+\sqrt{7}} \quad
 \qq{156}{3}=-\frac{\sqrt{2}}{3-\sqrt{7}} \quad
 \qq{124}{3}=-\frac{\sqrt{2}\,i\,}{3+\sqrt{7}}
\end{gather*}
and the map is given by
$$
\widetilde{h}(z)= \frac{-385 z^{14}+385 \sqrt{7} z^{12}-5159
z^{10}-2145 \sqrt{7} z^8+8085 \ z^6+795 \sqrt{7} z^4+315 z^2-3
\sqrt{7}}{-3 \sqrt{7} z^{15}-315 z^{13}+795 \ \sqrt{7} z^{11}-8085
z^9-2145 \sqrt{7} z^7+5159 z^5+385 \sqrt{7} z^3+385 z}.
$$
As a one-dimensional map $\widetilde{h}$ has $28=2\cdot15-2$
critical points.  Under restriction, the critical multiplicities of
the $21$ and $28$-points increase to four and three respectively.

\begin{figure}[ht]

\resizebox{\textwidth}{!}{\includegraphics{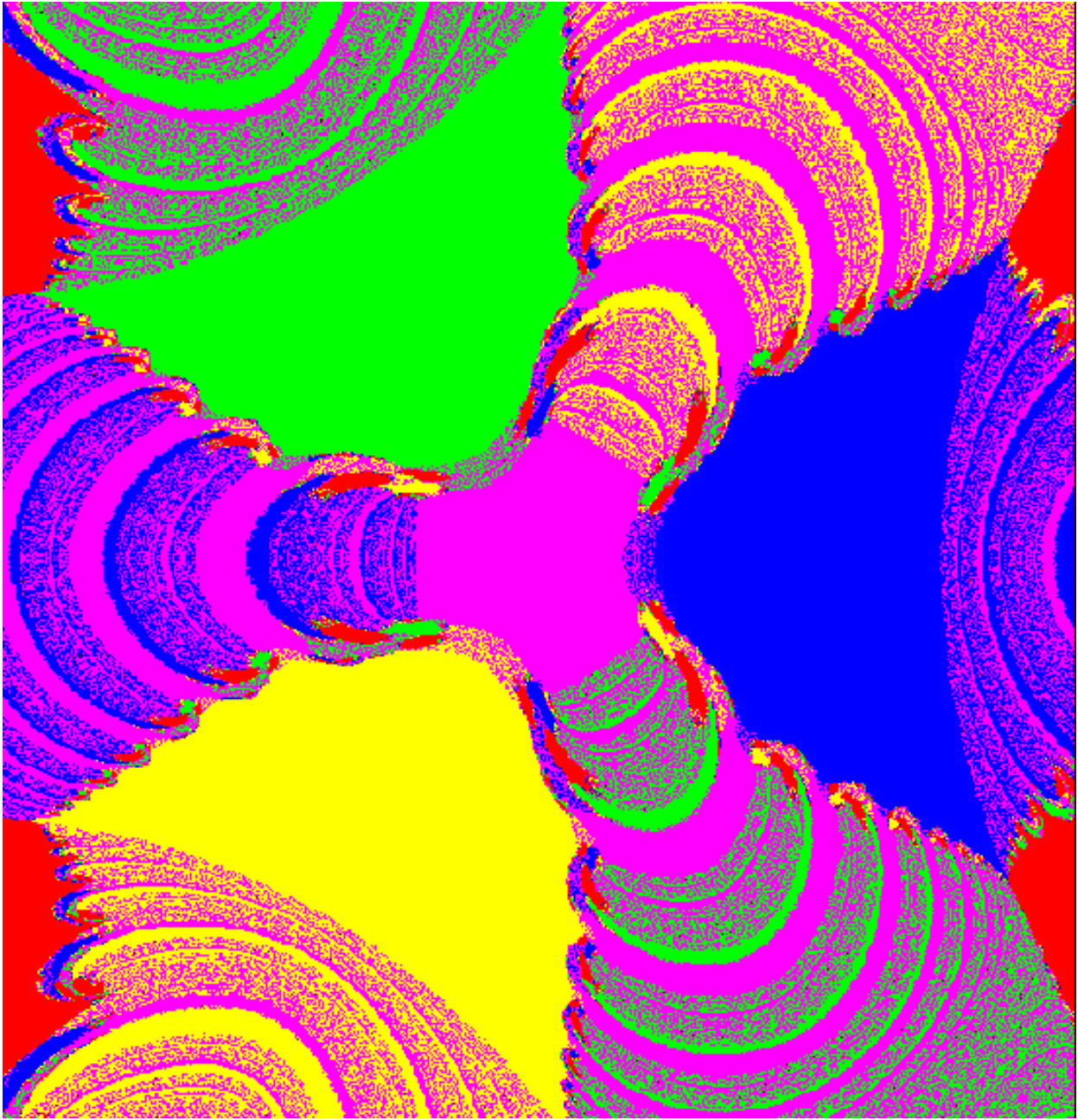}}

\caption{Basins of attraction on an \RP{2} mirror}

\label{fig:hRP2}

\end{figure}

Since its $y$-coefficients are real, $h$ preserves \RR{124}{7} and,
by symmetry, each reflection \RP{2}. Appearing in
Figure~\ref{fig:hRP2} is a basin plot for $h$ restricted to the
affine plane $\{y_3\neq 0\}$ of Figure~\ref{fig:quartic}. The
$y$-coordinates place the four $28$-points at
$$
 \qq{124}{7}=(0,0) \quad \qq{156}{3}=(1,0) \quad
 \qq{235}{6}=\biggl(-\frac{1}{2},\frac{\sqrt{3}}{2}\biggr) \quad
 \qq{346}{5}=\biggl(-\frac{1}{2},-\frac{\sqrt{3}}{2}\biggr)
$$
while the three $21$-points are on the line at infinity $\{y_3=0\}$:
$$
 \qq{37}{1}=[0,1,0] \quad
 \qq{67}{2}=[\sqrt{3},1,0] \quad
 \qq{57}{4}=[-\sqrt{3},1,0].
$$
Due to their location at infinity, the plot routine shows a single
basin for the $21$-points.

\clearpage
\section{Solving the heptic}

The treatment here follows that developed for the sextic.  A
detailed account appears in \cite{sextic}, \S\,4.

\subsection{Parameters for \K{168}}

To solve an equation, the fundamental task is to invert an
appropriate system of functions that shares the equation's symmetry.
The most familiar instance of this condition is the system of
elementary symmetric functions whose values are prescribed by a
polynomial's coefficients. In the case of a heptic with \K{168}
symmetry, the two degree-$0$ rational functions of
$$
K_1=\frac{\Phi(y)^2}{F(y)^3} \qquad K_2=\frac{\Psi(y)}{F(y) \Phi(y)}
$$
determine a single \K{168} orbit
$$\{K_1 F(y)^3 = \Phi(y)^2\}\cap  \{K_2 F(y) \Phi(y) = \Psi(y)\}$$
when values for $K_1$ and $K_2$ are chosen.  Note that assigning a
value to $K=(K_1,K_2)$ fails to produces every orbit; neither the
$24$-points (where $\{F=0\}\cap\{\Phi=0\}$) nor the $56$-points
(where $\{F=0\}\cap\{\Psi=0\}$) correspond to a value for $K$.
Finding the roots of a \K{168} heptic requires solving Klein's
\emph{form problem} which amounts to inverting these functions for a
specific $K$.

\subsection{Resolvent with \K{168} symmetry}

Taking the degree-$0$ rational functions
$$P_k(y)=\frac{C_k(y)^2}{F(y)} \quad (k=1\dots, 7)$$
to be the roots of a polynomial
$$
R_y(z)=\prod_{k=1}^7 (z-P_k(y))
 = \sum_{k=0}^7 b_k(y)\,z^{7-k},
$$
the $b_k$ are \Kl{2}-invariant---recall that the $P_k$ are permuted
by \Kl{2}. Accordingly, each coefficient is expressible as a
rational function in the basic invariants $F$, $\Phi$, and $\Psi$.
Converting to $K_1$ and $K_2$,
\small \begin{align*}
R_K(z)=\ & z^7- \frac{7}{48} (\theta +5) z^6 +
 \frac{7}{18432}(287\,\theta +507)z^5 \\
&+ \frac{7}{95551488}(1015\,\theta K_1 -420066\,\theta
-2253\,K_1-260010) z^4\\
&-\frac{49}{12230590464} (8460\,\theta K_1  - 1004571\,\theta
-16196\,K_1+66969) z^3 \\
&+\frac{7}{10567230160896} (9715653\,\theta K_1 + 69128\,\theta K_1
K_2 - 354295620\,\theta - 18385671\,K_1 \\
&- 195800\,K_1 K_2 +  228498732) z^2 \\
&-\frac{7}{36520347436056576} (5615155\,\theta K_1^2 -
15469601\,K_1^2 + 2374353432\,\theta K_1 + 43759152\,\theta K_2 K_1
\\
& - 132891408\,K_2 K_1 - 6013910664\,K_1 - 23865628752\,\theta +
32087745648)z
\\
&+\frac{1}{146081389744226304}(802165\,\theta K_1 K_2^2 -
2209943\,K_1 K_2^2 + 59819676\,\theta K_1 K_2 - 222370164\,K_1 K_2
\\
& + 964112436\,\theta K_1 - 5484159612\,K_1).
\end{align*} \normalsize
 An algorithm that harnesses the dynamics of the $15$-map
$h$ will compute a root of almost any member of this two-parameter
family of heptic resolvents.

\subsection{Self-parametrizing \K{168}}

Using \Kl{1} invariants and equivariants, define the coordinate
change
$$
u=T_y w=\Psi(y)\cdot y\,w_1 + \Phi(y)\cdot g_9(y)\,w_2 + F(y)\cdot
g_{11}(y)\, w_3
$$
that is parametrized by and has degree-$15$ in $y$.  By
construction, the map has an equivariance property in $y$:
$$T_{A y} = A\, T_y \quad \text{for all}\ A \in \Kl{1}.$$
Treating $u$ as \emph{reference} coordinates---with verbatim
replacement of $y$ by $u$, let $\K{}^u$ denote the $u$-expression of
\K{168} that acts on the $u$-space $\CP{2}_u$ so that the
$y$-parametrized family of \K{168} groups
$$\K{y}^w=T_y^{-1}\, \K{}^u\, T_y$$
acts on the $w$-space $\CP{2}_w$.

Each $\K{y}^w$ admits its system of invariants and equivariants.  In
the lowest-degree case,
$$F(u)=F(T_y w)$$
has degree $60=4\cdot 15$ in $y$ and degree $4$ in $w$. By
equivariance in $y$, the coefficient of each monomial in $w$ is
$\Kl{1}$ invariant. Expressing these coefficients in terms of
$F(y)$, $\Phi(y)$, and $\Psi(y)$ gives rise to a $K$-parametrized
form
$$F_K(w)=\frac{F(u)}{F(y)^{15}}.$$
Its explicit form appears at \cite{web}. The remaining invariant and
equivariant forms stem from $F_K$.

For clarity, let a variable subscript refer to differentiation.
Using transformation properties for the Hessian determinant and
defining
$$\Phi_K(w)=\alpha_\Phi \det{H_w(F_K(w))},$$
we obtain a $\K{y}^w$ invariant whose
$w$-degree is six:
\begin{align*}
\Phi(u)=&\ \alpha_\Phi\,\det{H_u(F(u))}=\
 \alpha_\Phi\,\det{H_u(F(T_y w))} \\
 =&\ \alpha_\Phi\,(\det{T_y})^{-2}\,\det{H_w(F(y)^{15} F_K(w))} \\
 =&\ (\det{T_y})^{-2}\,F(y)^{45}\,\Phi_K(w).
\end{align*}

Similar considerations yield forms of degrees $14$ and $21$:
\begin{align*}
\Psi(u)=&\ (\det{T_y})^{-6}\,F(y)^{120}\,\Psi_K(w) \\
X(u)=&\ (\det{T_y})^{-9}\,F(y)^{180}\,X_K(w).
\end{align*}

As for equivariant maps, a transformation formula for the gradient
operator produces
\begin{align*}
 \psi(u)=&\ \nabla_u F(u) \times \nabla_u \Phi(u) \\
 =&\ (\det{T_y})^{-3}\,F(y)^{60}\,T_y\,
  (\nabla_w F_K(w) \times \nabla_w \Phi_K(w)) \\
 =&\ (\det{T_y})^{-3}\,F(y)^{60}\,T_y\,\psi_K(w).
\end{align*}
With analogous definitions, maps in degrees $16$ and $18$ arise:
\begin{align*}
\phi(u)=&\ (\det{T_y})^{-7}\,F(y)^{135}\,T_y\,\phi_K(w) \\
f(u)=&\ (\det{T_y})^{-9}\,F(y)^{165}\,T_y\,f_K(w).
\end{align*}

We can now establish a $K$-parametrized version of the special
$15$-map $h$. As a preliminary issue, observe that by degree
counting and computation,
$$\det{T_y}=-2^43^27^{13}\, F \Phi \Psi X.$$
Now, define $\tau_K$ by
\begin{align*}
\frac{\det{T_y}^4}{F^{45}}=&\ 2^{16} 3^8 7^{52}\,\frac{(F\Phi\Psi
X)^4}{F^{45}} \\
=&\ 2^{16} 3^8 7^{52}\,\frac{X^4}{F^{21}}
\frac{\Phi^8}{F^{12}}\frac{\Psi^4}{F^8\Phi^4}\\
=&\ \tau_K K_1^4 K_2^4.
\end{align*}

By working directly in degree $36=15+21$, the unique map divisible
by $X_K(w)$ and doubly critical at $\{w\ |\ X_K(w)=0\}$ is given by
\small
\begin{align*}
h_{36}(u)=&\
 \bigl(-34992\,F(u)^7 + 6696\,F(u)^4\Phi(u)^2 + 1534\,F(u)\Phi(u)^4 +
 15912F(u)^2\Phi(u)\Psi(u) + 693\Psi(u)^2 \bigr) \psi(u) \\
 & + \bigl(5832\,F(u)^5 - 7290\,F(u)^2\Phi(u)^2 - 330\,\Phi(u)\Psi(u)
 \bigr) \phi(u) \\
 & + \bigl(5796\,F(u)^3\Phi(u) - 11\,\Phi(u)^3 + 216\,F(u)\Psi(u) \bigr)
 f(u)\\
=&\ F(y)^{300} (F(y)^{45} \tau_K)^{-\frac{15}{4}}\,T_y\Bigl(\bigl(
 -34992\,\tau_K^3 F_K(w)^7 + 6696\,\tau_K^2
 F_K(w)^4\Phi_K(w)^2  \\
 & + 1534\,\tau_K F_K(w)\Phi_K(w)^4 + 15912\tau_K
 F_K(w)^2\Phi_K(w)\Psi_K(w) + 693\Psi(w)^2 \bigr)
 \psi_K(w) \\
 & + \bigl(5832\,\tau_K^2 F_K(w)^5 - 7290\,\tau_K
 F_K(w)^2\Phi_K(w)^2 - 330\,\Phi_K(w)\Psi_K(w)
 \bigr) \phi_K(w) \\
 & + \bigl(5796\,\tau_K F_K(w)^3\Phi_K(w) - 11\,\Phi_K(w)^3 +
 216\,F_K(w)\Psi_K(w) \bigr) f_K(w) \Bigr).
\end{align*}
\normalsize
For a map on $\CP{2}_w$, take
\small
\begin{align*}
h_{36_K}(w)=&\ \bigl(
 -34992\,\tau_K^3 F_K(w)^7 + 6696\,\tau_K^2
 F_K(w)^4\Phi_K(w)^2 + 1534\,\tau_K F_K(w)\Phi_K(w)^4 \\
 &  + 15912\tau_K F_K(w)^2\Phi_K(w)\Psi_K(w) + 693\Psi(w)^2 \bigr)
 \psi_K(w) \\
 & + \bigl(5832\,\tau_K^2 F_K(w)^5 - 7290\,\tau_K
 F_K(w)^2\Phi_K(w)^2 - 330\,\Phi_K(w)\Psi_K(w)
 \bigr) \phi_K(w) \\
 & + \bigl(5796\,\tau_K F_K(w)^3\Phi_K(w) - 11\,\Phi_K(w)^3 +
 216\,F_K(w)\Psi_K(w) \bigr) f_K(w)
\end{align*}
\normalsize
so that the $15$-map is
$$h_K(w)= \frac{h_{36_K}(w)}{X_K(w)}.$$
The final task is to harness $h_K$'s dynamics to a root-finding
procedure.

\subsection{From dynamics to a root}

Direct computation yields (octahedral) \Oct{\ell} invariants
$A_\ell$ in degree-$10$ such that 1) the $A_\ell$ are simply
permuted by $\Kl{1}$ and 2) at a $21$ or $28$-point, all but one of
the $A_\ell$ vanish.  Furthermore,
$$
 B_\ell(y)=
 \frac{2401(57+11\sqrt{7}\,i)}{18432}\,\frac{A_\ell(y)}{F(y)\Phi(y)}
$$
is the degree-$0$ rational function normalized so that it takes on
the value $1$ at those $21$ or $28$-points where it fails to vanish.
Explicitly,
\begin{align*}
B_\ell(\qq{ac}{b})=\begin{cases}
0&\ell\neq b\\
1&\ell=b
\end{cases} \qquad
B_\ell(\qq{abc}{d})=\begin{cases}
0&\ell\neq d\\
1&\ell=d
\end{cases}.
\end{align*}

A choice of values $K=(K_1,K_2)$ gives a resolvent $R_K$ whose roots
are $P_\ell(y)$. Define the $\K{}^y$-invariant rational function
$$J_y(w)= \sum_{\ell=1}^7 B_\ell(T_y w) P_\ell(y).$$
Being degree-$0$ in $y$, $J_y(w)$ converts to $J_K(w)$ whose
$w$-coefficients are expressions in $K$. Given a $21$ or $28$-point
$W$ in $\CP{2}_w$, there is a $21$ or $28$-point $U=T_y W$ in
$\CP{2}_u$.  When evaluated at $W$, all but one value $B_\ell(T_y
W)=B_\ell(U)$ vanishes---when $\ell=m$, say.  Thus,
$$J_K(W)=B_m(U)P_m(y)=P_m(y)$$
is a root of $R_K$.

Hence, a root of $R_K$ results from finding a $21$ or $28$-point in
$\CP{2}_w$.  The dynamics of $h_K$ provides a means to this end:
from a symmetry-breaking random start $w_0$, iteration of $h_K$
locates a $21$ or $28$-point under the $\K{y}^w$ action.

\emph{Mathematica} files that implement the heptic-solving algorithm
are available at \cite{web}.

\subsection{Obtaining a second root}

In order to locate a $21$ or $28$-point in $\CP{2}_w$, the heptic
algorithm above relies on their attracting nature under $h_K$.
However, as a trajectory $\{h_K^\ell(w_0)\}$ nears a $21$ or
$28$-point it does so \emph{along} one of the attracting $21$-lines.
We can use this piece of line-specific information to compute a
second root of the resolvent $R_K$.  What we need is a means of
finding a second $21$ or $28$-point $V$ such that $J_K(V)$ is a root
different from the first. The special behavior of $\psi$ gives us a
tool for the task.

Since $\psi=\nabla F \times \nabla \Phi,$ the map's geometric effect
on a point $a$ produces a second point $\psi(a)$ as the intersection
of two lines:
$$\psi(a)=\{\nabla F(a)^T y =0\} \cap \{\nabla \Phi(a)^T y =0\}.$$
When $a$ lies on a $21$-line \LL{ac}{b}, the lines that determine
$\psi(a)$ each pass through \qq{ac}{b}, a concurrence that reflects
the \K{168} duality between points on a $21$-line and the pencil of
lines through its companion point. Thus,
$\psi(\LL{ac}{b})=\qq{ac}{b}$. This ``blowing-down" of the
$21$-lines requires their intersections---the $21$ and
$28$-points---to blow up, that is, for any lift $\widehat{\psi}$ and
$\widehat{q}$ of $\psi$ and a $21$ or $28$-point $q$,
$\widehat{\psi}(\widehat{q})=0$.

By watching $\psi_K(h_K^\ell(w_0))$ for fixed point behavior, we can
tell when $h_K^\ell(w_0)$ has found a $21$-line.  Then
$J_K(\psi_K(h_K^\ell(w_0)))$ approximates a root of $R_K$.
Unfortunately, this might give the same root as $J_K(W)$ where
$W=\lim_{\ell \rightarrow \infty}h_K^\ell(w_0)$.  For instance, let
the $21$-line be \LL{13}{7} which contains the $21$-points
\qq{26}{7} and \qq{45}{7}. When evaluated at either of these points,
the root selector $J_K$ gives the same result as it does at
\qq{13}{7}, the image of the blown-down line.  So, if the iteration
finds one of the other two $21$-points or one of the four
$28$-points on \LL{13}{7}, blowing down the line provides for a
second root.

In practice, this tactic requires some vigilance.  Applying $\psi_K$
is effective only after the trajectory reaches the vicinity of a
$21$-line but \emph{before} it closes in on a $21$ or $28$-point. By
allowing $h_K^\ell(w_0)$ near a $21$ or $28$-point, the blow-up of
the points onto their companion $21$ or $28$-line overwhelms the
blow-down of the $21$-line.



\begin{thebibliography}{99}

\bibitem[C1]{sextic} S.\ Crass. \emph{Solving the sextic by
iteration: A study in complex geometry and dynamics.} Experiment.\
Math.\ 8 (1999) No.\ 3, 209-240. Preprint at
\texttt{arxiv.org/abs/math.DS/9903111}

\bibitem[C2]{critFin} S.\ Crass. \emph{A family of critically finite
maps with symmetry.}  Publ.\ Mat.\ 49 (2005), No.\ 1, 127--157.
Preprint at \texttt{arxiv.org/abs/math.DS/0307057}

\bibitem[C3]{web} S.\ Crass.
\texttt{www.csulb.edu/$\sim$scrass/Math}

\bibitem[F]{fricke} R.\ Fricke. \emph{Lehrbuch der Algebra} 2.
Vieweg, 1926.

\bibitem[K]{klein} F.\ Klein.  \emph{\"{U}ber die Transformation
siebenter Ordnung der elliptischen Funktionen.}  Mathematische
Annalen 14 (1879), 428-471.  Translation by S.\ Levy appears as
\emph{On the Order-Seven Transformation of Elliptic Functions} in
\cite{8fold}.

\bibitem[L]{8fold} S.\ Levy, ed.  \emph{The Eightfold Way} MSRI
Publ.\ Vol.\ 35 (1998).

\bibitem[ST]{st} G.\ Shephard and T.\ Todd.
\emph{Finite unitary reflection groups.} Canad.\ J.\ Math.\ 6
(1954), 274-304.

\bibitem[S]{shepler} A.\ Shepler. \emph{Semi-invariants of finite
reflection groups}. J.\ Algebra 220 (1999), 314-326.

\end{thebibliography}
\end{document}